%
%
%
%
\hsize=5in
\baselineskip=12pt
\vsize=20cm
\parindent=10pt
\pretolerance=40
\predisplaypenalty=0
\displaywidowpenalty=0
\finalhyphendemerits=0
\hfuzz=2pt
\frenchspacing
\footline={\ifnum\pageno=1\else\hfil\tenrm\number\pageno\hfil\fi}
%
%
\input amssym.def
\def\titlefonts{\baselineskip=1.44\baselineskip
	\font\titlef=cmbx12
	\titlef
	}
\font\ninerm=cmr9
\font\ninebf=cmbx9
\font\ninei=cmmi9
\skewchar\ninei='177
\font\ninesy=cmsy9
\skewchar\ninesy='60
\font\nineit=cmti9
\def\reffonts{\baselineskip=0.9\baselineskip
	\textfont0=\ninerm
	\def\rm{\fam0\ninerm}%
	\textfont1=\ninei
	\textfont2=\ninesy
	\textfont\bffam=\ninebf
	\def\bf{\fam\bffam\ninebf}%
	\def\it{\nineit}%
	}
%
%
\def\frontmatter{\vbox{}\vskip1cm\bgroup
	\leftskip=0pt plus1fil\rightskip=0pt plus1fil
	\parindent=0pt
	\parfillskip=0pt
	\pretolerance=10000
	}
\def\endfrontmatter{\egroup\bigskip}
\def\title#1{{\titlefonts#1\par}}
\def\author#1{\bigskip#1\par}
\def\address#1{\medskip{\reffonts\it#1}}
\def\email#1{\smallskip{\reffonts{\it E-mail: }\rm#1}}
\def\thanks#1{\footnote{}{\reffonts\rm\noindent#1\hfil}}
\def\section#1\par{\ifdim\lastskip<\bigskipamount\removelastskip\fi
	\penalty-250\bigskip
	\vbox{\leftskip=0pt plus1fil\rightskip=0pt plus1fil
	\parindent=0pt
	\parfillskip=0pt
  \pretolerance=10000{\bf#1}}\nobreak\medskip
	}
\def\proclaim#1. {\medbreak\bgroup{\noindent\bf#1.}\ \it}
\def\endproclaim{\egroup
	\ifdim\lastskip<\medskipamount\removelastskip\medskip\fi}
\newdimen\itemsize
\def\setitemsize#1 {{\setbox0\hbox{#1\ }
	\global\itemsize=\wd0}}
\def\item#1 #2\par{\ifdim\lastskip<\smallskipamount\removelastskip\smallskip\fi
	{\leftskip=\itemsize
	\noindent\hskip-\leftskip
	\hbox to\leftskip{\hfil\rm#1\ }#2\par}\smallskip}
\def\Proof#1. {\ifdim\lastskip<\medskipamount\removelastskip\medskip\fi
	{\noindent\it Proof\if\space#1\space\else\ \fi#1.}\ }
\def\endproof{\hfill\hbox{}\quad\hbox{}\hfill\llap{$\square$}\medskip}
\def\Remark #1.{\ifdim\lastskip<\medskipamount\removelastskip\medskip\fi
        {\noindent\bf Remark #1.}}
\def\endremark{\medskip}
\def\emph#1{{\it #1}\/}
%
%
\newcount\citation
\newtoks\citetoks
\def\citedef#1\endcitedef{\citetoks={#1\endcitedef}}
\def\endcitedef#1\endcitedef{}
\def\citenum#1{\citation=0\def\curcite{#1}%
	\expandafter\checkendcite\the\citetoks}
\def\checkendcite#1{\ifx\endcitedef#1?\else
	\expandafter\lookcite\expandafter#1\fi}
\def\lookcite#1 {\advance\citation by1\def\auxcite{#1}%
	\ifx\auxcite\curcite\the\citation\expandafter\endcitedef\else
	\expandafter\checkendcite\fi}
\def\cite#1{\makecite#1,\cite}
\def\makecite#1,#2{[\citenum{#1}\ifx\cite#2]\else\expandafter\clearcite\expandafter#2\fi}
\def\clearcite#1,\cite{, #1]}
%
%
\def\references{\section References\par
	\bgroup
	\parindent=0pt
	\reffonts
	\rm
	\frenchspacing
	\setbox0\hbox{99. }\leftskip=\wd0
	}
\def\endreferences{\egroup}
\newtoks\authtoks
\newif\iffirstauth
\def\checkendauth#1{\ifx\auth#1%
    \iffirstauth\the\authtoks
    \else{} and \the\authtoks\fi,%
  \else\iffirstauth\the\authtoks\firstauthfalse
    \else, \the\authtoks\fi
    \expandafter\nextauth\expandafter#1\fi
	}
\def\nextauth#1,#2;{\authtoks={#1 #2}\checkendauth}
\def\auth#1{\nextauth#1;\auth}
\newif\ifinbook
\newif\ifbookref
\def\nextref#1 {\par\hskip-\leftskip
	\hbox to\leftskip{\hfil\citenum{#1}.\ }%
	\initnextref}
\def\initnextref{\bookreffalse\inbookfalse\firstauthtrue\ignorespaces}
\def\paper#1{{\it#1},}
\def\InBook#1{\inbooktrue in ``#1",}
\def\book#1{\bookreftrue{\it#1},}
\def\journal#1{#1\ifinbook,\fi}
\def\BkSer#1{#1,}
\def\Vol#1{{\bf#1}}
\def\BkVol#1{Vol. #1,}
\def\publisher#1{#1,}
\def\Year#1{\ifbookref #1.\else\ifinbook #1,\else(#1)\fi\fi}
\def\Pages#1{\makepages#1.}
\long\def\makepages#1-#2.#3{\ifinbook pp. \fi#1--#2\ifx\par#3.\fi#3}
\def\inRus{{ \rm(in Russian)}}
%
%
\newsymbol\square 1003
\newsymbol\varnothing 203F
\let\Rar\Rightarrow
\let\Lrar\Leftrightarrow
\let\hrar\hookrightarrow
\let\lhu\leftharpoonup
\let\rhu\rightharpoonup
\let\ot\otimes
\let\sbs\subset
\let\sps\supset
\def\cop{^{\rm cop}}
\def\End{\mathop{\rm End}\nolimits}
\def\GKdim{\mathop{\rm GKdim}\nolimits}
\def\Hom{\mathop{\rm Hom}\nolimits}
\def\id{{\rm id}}
\def\Ker{\mathop{\rm Ker}}
\def\lann{\mathop{\rm lann}\nolimits}
\def\limdir{\mathop{\vtop{\offinterlineskip\halign{##\hskip0pt\cr\rm lim\cr
	\noalign{\vskip1pt}
	$\scriptstyle\mathord-\mskip-10mu plus1fil
	\mathord-\mskip-10mu plus1fil
	\mathord\rightarrow$\cr}}}}
\def\lng{\mathop{\rm length}\nolimits}
\def\op{^{\rm op}}
\def\rann{\mathop{\rm rann}\nolimits}
\def\Spec{\mathop{\rm Spec}}
\def\trdeg{\mathop{\rm tr\,deg}}
\def\mapr#1{{}\mathrel{\smash{\mathop{\longrightarrow}\limits^{#1}}}{}}
\let\al\alpha
\let\be\beta
\let\ep\varepsilon
\let\la\lambda
\let\th\theta
\let\De\Delta
\let\Om\Omega
\def\0{_{(0)}}
\def\1{_{(1)}}
\def\2{_{(2)}}
\def\3{_{(3)}}
\def\C{{\cal C}}
\def\F{{\cal F}}
\def\G{{\cal G}}
\def\Hd{H^\circ}
\def\I{{\cal I}}
\def\L{{\cal L}}
\def\T{{\cal T}}
\citedef
An92
Br84
Br07
Chir
Coh85
Dem-G
Doi93
Goo
Goo-W
Gro
Her-S64
Kr-L
Mal43
Ma91
Ma-W94
Mc-R
Mol75
Mo
Mo-Sch95
Mu-Sch99
Nich-Z89
Pro
Rad77
Row
Scha92
Scha00
Schn92
Schn93
Sk06
Sk07
Sk08
Sk10
Sk99
Sk-Oy06
St
Tak79
Wu-Zh03
\endcitedef

\frontmatter

\title{Flatness over PI coideal subalgebras}
\author{Serge Skryabin}
\address{Institute of Mathematics and Mechanics,
Kazan Federal University, Russia\break}
\email{Serge.Skryabin@kpfu.ru}

\endfrontmatter

\section
Introduction

Construction of algebraic homogeneous spaces for group schemes of finite type 
over a field is linked closely with the flatness of orbit morphisms \cite{Dem-G}. 
This motivated the study of the flatness property for Hopf algebras. At a 
certain time it seemed plausible that an arbitrary Hopf algebra over a field 
might be faithfully flat over any Hopf subalgebra (see \cite{Mo, Question 
3.5.4}). This is not true, as was shown by Schauenburg \cite{Scha00}. Still, 
faithful flatness holds under additional assumptions, although there is no any 
single result which would unify the known cases.

To interpret homogeneous spaces one needs not only Hopf subalgebras but the 
larger class of right coideal subalgebras. Over coideal subalgebras one can 
expect flatness in favourable cases but not necessarily faithful flatness. In 
the author's preceding article \cite{Sk99} it was proved that a residually 
finite-dimensional Noetherian Hopf algebra is a flat left module over any 
right Noetherian right coideal subalgebra.

An obvious deficiency of that result is that nothing is said about flatness 
when the Noetherian condition on the subalgebra is dropped. Even a finitely 
generated commutative Hopf algebra may contain non-Noetherian coideal 
subalgebras. Such a possibility stems from the existence of algebraic groups 
admitting a finite-dimensional linear representation such that the respective 
algebra of invariant polynomial functions is not finitely generated (see 
\cite{Gro}). 

In the new paper we deal with right coideal subalgebras satisfying a polynomial 
identity. Polynomial identity algebras, PI algebras for short, constitute 
an important class which includes, in particular, all commutative algebras and 
all algebras module-finite over some of their commutative subalgebras.

We fix a field $k$ which is assumed to be the base field for all algebras. 
Recall that an algebra is said to be \emph{residually finite dimensional} if 
its ideals of finite codimension have zero intersection \cite{Mo}.

\proclaim
Theorem 0.1.
A residually finite-dimensional Noetherian Hopf algebra $H$ is left and right 
flat over any PI right coideal subalgebra and is left and right faithfully flat 
over any PI Hopf subalgebra.
\endproclaim

By a result of Schneider \cite{Schn92, Cor. 1.8}, for any Hopf algebra with 
bijective antipode faithful flatness over a Hopf subalgebra is equivalent to 
projectivity (see also Masuoka and Wigner \cite{Ma-W94}). This applies in the 
situation of Theorem 0.1, so in the case of Hopf subalgebras we have an even 
stronger conclusion. If $H$ is additionally assumed to be itself PI, then the 
polynomial identity is inherited by all subalgebras of $H$, which yields

\proclaim
Corollary 0.2.
A residually finite-dimensional Noetherian PI Hopf algebra $H$ is left and 
right flat over any right coideal subalgebra and is a projective generator as 
either left or right module over any Hopf subalgebra.
\endproclaim

This conclusion extends immediately to the case when the Hopf algebra $H$ is 
the union of a directed family of residually finite-dimensional Noetherian PI 
Hopf subalgebras. In such a form our result covers two classes of Hopf algebras 
which have been previously known to be flat over all coideal subalgebras: the 
commutative Hopf algebras \cite{Ma-W94} and the finite-dimensional ones \cite{Sk07}.

Faithful flatness and projectivity over all Hopf subalgebras have been known 
in the two cases just mentioned \cite{Tak79}, \cite{Nich-Z89}, and in two 
other cases rather different by the arguments employed: the Hopf algebras with 
cocommutative coradicals \cite{Ma91} (including pointed Hopf algebras 
\cite{Rad77}) and the cosemisimple ones \cite{Chir}. Partial results on 
faithful flatness over PI Hopf subalgebras of some type can be found in 
\cite{Schn93}, \cite{Sk08}, \cite{Wu-Zh03}.

Over coideal subalgebras the flatness property turns out to be more elusive. 
There are several results in which faithful flatness over some coideal 
subalgebras has been proved, but this requires usually more severe 
restrictions than in the case of Hopf subalgebras (see \cite{Ma91}, 
\cite{Mu-Sch99}, \cite{Sk08}).

Theorem 0.1 is another instance of the approach developed in \cite{Sk10}. In 
that paper it was shown that a Hopf algebra $H$ is left flat over a right 
coideal subalgebra $A$ whenever both $H$ and $A$ have right Artinian classical 
right quotient rings (the Ore rings of fractions). For a residually 
finite-dimensional Noetherian Hopf algebra $H$ and its right Noetherian right 
coideal subalgebras the existence of such quotient rings has been established 
in \cite{Sk99}. Here we deal with PI coideal subalgebras. By passage to the 
dual Hopf algebra the whole problem is reformulated in the context of 
\hbox{$H$-prime} $H$-module algebras where $H$ is now an arbitrary Hopf 
algebra over the base field $k$.

As in \cite{Sk99}, we aim first to construct a generalized quotient ring $Q(A)$ 
and then to verify that $Q(A)$ is the Ore localization of $A$. However, the 
right Gabriel topology used in \cite{Sk99} is well suited only when $A$ is 
right Noetherian. In the present paper we use a different Gabriel topology 
defined in terms of Gelfand-Kirillov dimension. It is characterized by 
means of the corresponding class of torsion modules.

Another essential difference with \cite{Sk99} is that we have to work with a 
symmetric quotient ring $Q(A)$ defined with respect to a pair of right and left 
Gabriel topologies. It can be compared to the symmetric Martindale quotient ring, 
but in our case two filters of one-sided ideals are involved in the construction. 
Although extra work is needed to verify that $Q(A)$ is semiprimary, once this 
is done we gain considerable advantage since the situation becomes completely 
left-right symmetric, unlike what we had in \cite{Sk99}.

The first two sections of the paper provide a purely ring-theoretic background. 
For a finitely generated PI algebra $A$ we introduce certain torsion classes 
$\T_r$ and $\T_l$ of right and left $A$-modules and prove several lemmas aimed 
at recognition of torsion modules. The Gabriel topologies $\G_r$ and $\G_l$ 
correspond to $\T_r$ and $\T_l$.

Section 2 is concerned with the $(\G_l,\G_r)$-symmetric quotient ring $Q(A)$. 
Proposition 2.9 gives a set of conditions under which $Q(A)$ is shown to be an 
Artinian classical quotient ring of $A$. In such a form this result is probably 
not very useful since the ring-theoretic assumptions are rather demanding.

However, all necessary conditions can be verified for a finitely generated 
$H$-prime $H$-module PI algebra $A$ satisfying the ACC (ascending chain 
condition) on right and left annihilators when the action of $H$ on $A$ is 
locally finite. This is done in section 3, and the final result is presented 
in Theorem 3.9. Its extension to the \hbox{$H$-semiprime} case is given in 
Theorem 3.12. In the $H$-prime case it is shown that the algebra $A$ behaves 
particularly well with respect to the Gelfand-Kirillov dimension.

In section 4 we apply Theorem 3.9 to coideal subalgebras. Theorem 4.1 is the 
main result of the paper which includes Theorem 0.1 as a special case. With 
the assumption that $H$ has an Artinian classical quotient ring the 
Noetherianness of $H$ is not needed.

Theorem 4.5 is another special case of Theorem 4.1 stated for a PI Hopf 
algebra $H$ which is \emph{affine}, i.e., finitely generated as an ordinary 
algebra. Here the Noetherian condition on $H$ is relaxed to the ACC on 
annihilators. We are led to suggest that the conclusion might hold 
even in larger generality:

\proclaim
Conjecture 0.3.
Every affine PI Hopf algebra $H$ is flat over all right coideal subalgebras 
and is faithfully flat over all Hopf subalgebras.
\endproclaim

Among all Hopf algebras those satisfying a polynomial identity are more 
manageable, but even this class of Hopf algebras is not understood 
sufficiently well. Let us look at the following assertions:

\setitemsize(a)
\item(a)
every affine PI Hopf algebra is Noetherian,

\item(b)
every affine PI Hopf algebra is \emph{representable} in the sense that it 
embeds in a finite-dimensional algebra over an extension field of the ground 
field,

\item(c)
every affine PI Hopf algebra is residually finite dimensional.

\smallskip
It is not known whether any of them is true or not. As to (a), this was asked 
by Brown \cite{Br07, Question C}. The implications (a)$\,\Rar\,$(b)$\,\Rar\,$(c) 
follow from results of Anan'in \cite{An92} and Malcev \cite{Mal43}. Knowledge 
of (b) would be sufficient to confirm Conjecture 0.3.

In the direction opposite to (a) Wu and Zhang asked whether every Noetherian 
PI Hopf algebra is affine \cite{Wu-Zh03, Question 5.1}. Proposition 4.8 and 
Corollary 4.9 in our paper are partial results on that question deduced 
quickly from the main result.

\section
1. The torsion determined by the Gelfand-Kirillov dimension

Let $A$ be a finitely generated PI algebra. By an algebra we mean an 
associative unital algebra over a base field $k$. We will use standard facts 
from PI theory discussed thoroughly in \cite{Mc-R}, \cite{Pro}, \cite{Row}. 
Concerning the definition and properties of Gelfand-Kirillov dimension we 
refer to \cite{Kr-L}, \cite{Mc-R}, \cite{Row}.

The algebra $A$ has finite Gelfand-Kirillov dimension $\GKdim A$ \cite{Row, 
Th. 6.3.25}. For an $A$-module $M$ we denote its Gelfand-Kirillov dimension 
by $\GKdim M$ or by $\GKdim_AM$ when there is a need to specify the algebra.
We have
$$
\GKdim M\le\GKdim A. 
$$
If $\GKdim M'=\GKdim M$ for all nonzero submodules $M'\sbs M$, then $M$ is 
called \emph{GK-pure} (or \emph{GK-homogeneous}). In particular, $A$ is 
\emph{right GK-pure} if
$$
\GKdim I=\GKdim A
$$
for all nonzero right ideals $I$ of $A$. It is known that each prime right 
Goldie algebra is right GK-pure \cite{Kr-L, Lemma 5.12}.

If $I$ is an ideal of $A$ then the Gelfand-Kirillov dimension of $A/I$ 
as an algebra is the same as that of $A/I$ as either right or left $A$-module. Put
$$
\openup1\jot
\eqalign{
d(A)&{}=\max\{\GKdim A/P\mid P\in\Spec A\},\cr
\Om(A)&{}=\{P\in\Spec A\mid\GKdim A/P=d(A)\}
}
$$
where $\Spec A$ is the set of prime ideals of $A$.

If $P'\sbs P$ is a proper inclusion between two prime ideals of $A$, then 
$\GKdim A/P<\GKdim A/P'$ by \cite{Kr-L, Prop. 3.15} since the ideal $P/P'$ of 
the prime factor algebra $A/P'$ contains a regular element (nonzerodivisor) of 
the latter. Therefore all ideals in the set $\Om(A)$ are minimal primes of $A$. 
It is known that $A$ has finitely many minimal prime ideals \cite{Pro, Ch. 5, 
Cor. 2.2, 2.4}. Hence $\Om(A)$ is a finite set.

The prime radical $N$ of $A$ is the intersection of all prime ideals of 
$A$. By \cite{Kr-L, Cor. 3.3} we have
$$
d(A)=\GKdim A/N\le\GKdim A.
$$

Denote by $\T_r$ the class of all right $A$-modules $M$ which admit a finite 
chain of submodules $M=M_0\sps M_1\sps\cdots\sps M_n=0$ with all factors 
$M_{i-1}/M_i$ of Gelfand-Kirillov dimension less than $d(A)$.

A ring $Q$ is said to be a \emph{classical right quotient ring} (resp., 
\emph{left quotient ring}) of a ring $R$ if $R$ embeds in $Q$ as a subring, all 
regular elements (nonzerodivisors) of $R$ are invertible in $Q$ and each element 
of $Q$ can be written as $as^{-1}$ (resp., as $s^{-1}a$) for some $a,s\in R$ with 
$s$ being regular. If the right and left conditions in this definition are 
satisfied simultaneously, then $Q$ is a \emph{classical quotient ring} of $R$.

Rings admitting a semisimple Artinian classical one-sided quotient ring were 
characterized by Goldie. If $P\in\Spec A$, then the prime factor ring $A/P$ is 
right and left Goldie by Posner's Theorem. Therefore $A/P$ has a simple 
Artinian classical quotient ring $Q(A/P)$.

\proclaim
Lemma 1.1.
For a right $A$-module $M$ annihilated by a prime ideal $P\in\Om(A)$ 
the following conditions are equivalent{\rm:}

\item(a)
\ $M\in\T_r\mskip1mu,$

\item(b)
\ $\GKdim M<d(A)\mskip1mu,$

\item(c)
\ $M\ot_AQ(A/P)=0\mskip1mu$.

\endproclaim

\Proof.
If $M\ot_AQ(A/P)=0$, then each element of $M$ is annihilated by a regular 
element of $A/P$, i.e., $M$ is a torsion right $A/P$-module in the classical 
sense. In this case $\GKdim M<\GKdim A/P$ by \cite{Mc-R, 8.3.6} since 
the ring $A/P$ is right Goldie. Thus (c)$\Rar$(b). The implication (b)$\Rar$(a) 
is obvious.

Suppose that $M\ot_AQ(A/P)\ne0$. Since each $Q(A/P)$-module is semisimple, and 
each simple right $Q(A/P)$-module is isomorphic to a right ideal of $Q(A/P)$, 
we have $\Hom_A\bigl(M,\,Q(A/P)\bigr)\ne0$. In other words, some nonzero 
epimorphic image $M'$ of $M$ is isomorphic to a right $A$-submodule of 
$Q(A/P)$. But each nonzero $A$-submodule of $Q(A/P)$ has nonzero intersection 
with $A/P$. Therefore there is a nonzero right ideal $I$ of $A/P$ isomorphic 
to a submodule of $M'$. Recalling that prime right Goldie algebras are 
right GK-pure \cite{Kr-L, Lemma 5.12}, we get
$$
\GKdim M\ge\GKdim M'\ge\GKdim I=\GKdim A/P=d(A).
$$

The preceding argument shows that (b)$\Rar$(c). If
$$
M=M_0\sps M_1\sps\cdots\sps M_n=0
$$
is a chain of submodules such that $\GKdim M_{i-1}/M_i<d(A)$ for each 
$i=1,\ldots,n$, then $M_{i-1}/M_i\ot_AQ(A/P)=0$ for each $i$, which entails 
$M\ot_AQ(A/P)=0$. Thus (a)$\Rar$(c) too.  
\endproof

If $M$ is annihilated by a prime ideal $P$ of $A$ such that $P\notin\Om(A)$, then
$$
\GKdim M\le\GKdim A/P<d(A).
$$
In this case $\,M\in\T_r\mskip1mu$.

\proclaim
Lemma 1.2.
The class $\T_r$ is closed under submodules, factor modules, extensions and 
coproducts. In other words{\rm,} $\T_r$ is the torsion class of a hereditary 
torsion theory\/ {\rm(see \cite{St, Ch. VI, \S3})}.
\endproclaim

\Proof.
Given an exact sequence of right $A$-modules $0\to M'\to M\to M''\to0$, it is 
clear from the definition that $M\in\T_r$ whenever $M'\in\T_r$ and $M''\in\T_r$. 
Conversely, any chain of submodules of $M$ with factors of Gelfand-Kirillov 
dimension less than $d(A)$ induces similar chains in $M'$ and $M''$ since the 
Gelfand-Kirillov dimensions of submodules and factor modules of any module do 
not exceed the dimension of the module itself. Thus $M\in\T_r$ if and only 
if $M'\in\T_r$ and $M''\in\T_r$.

It is known that the prime radical $N$ of $A$ is nilpotent \cite{Br84}. 
Since the product of the minimal prime ideals of $A$ is contained in $N$, 
there exists a finite sequence of prime ideals $P_1,\ldots,P_n$ such that 
$P_1\cdots P_n=0$. Setting $K_0=A$ and $K_i=K_{i-1}P_i$ for $i>0$, we obtain a 
descending chain of ideals ending at $K_n=0$. For an arbitrary right 
$A$-module $M$ we get a chain of submodules
$$
M=M_0\sps M_1\sps\cdots\sps M_n=0\quad\hbox{with $M_i=MK_i$}.
$$
By the already established properties $M\in\T_r$ if and only if 
$M_{i-1}/M_i\in\T_r$ for each $i=1,\ldots,n$. Note that $M_{i-1}/M_i$ is 
annihilated by $P_i$. If $P_i\notin\Om(A)$, then $\GKdim M_{i-1}/M_i<d(A)$. 
Taking into account Lemma 1.1, we deduce that $M\in\T_r$ if and only if 
$\,\GKdim M_{i-1}/M_i<d(A)\,$ for each $i=1,\ldots,n$.

Suppose that $M$ is the direct sum of a collection of right $A$-modules 
$\{T_\al\}\sbs\T_r$. By the previous observation 
$\GKdim T_\al K_{i-1}/T_\al K_i<d(A)$ for each $\al$ and $i$. Since
$$
MK_{i-1}/MK_i\cong\bigoplus_\al\,T_\al K_{i-1}/T_\al K_i\,,
$$
we get $\GKdim MK_{i-1}/MK_i<d(A)$ for each $i$, whence $M\in\T_r$. Thus for any 
set of modules in $\T_r$ their direct sum also lies in $\T_r$.
\endproof

The modules in $\T_r$ can be characterized by means of a variation of the notion 
of reduced rank (see \cite{Goo-W, Ch. 11} and \cite{Mc-R, Ch. 4}). Let $P$ be 
any minimal prime ideal of $A$. Recall that a \emph{Loewy series} for a right 
$A$-module $M$ is any finite chain of submodules $M=M_0\sps M_1\sps\cdots\sps 
M_n=0$ with all factors annihilated by the prime radical $N$ of $A$. Define the 
\emph{reduced rank} $\rho_P(M)$ of $M$ at $P$ as follows. Choose any Loewy 
series for $M$ and put
$$
\rho_P(M)=\sum_{i=1}^n\,\lng M_{i-1}/M_i\ot_AQ(A/P)
$$
where the summands are the lengths of right $Q(A/P)$-modules. It is assumed 
that $\rho_P(M)=+\infty$ when at least one of the $Q(A/P)$-modules 
$M_{i-1}/M_i\ot_AQ(A/P)$ is not finitely generated.

The number $\rho_P(M)$ does not depend on the choice of a Loewy series. If 
the whole $M$ is annihilated by $N$, this follows from the equality
$$
\sum_{i=1}^n\,\lng M_{i-1}/M_i\ot_AQ(A/P)=\lng M\ot_AQ(A/P)
$$
which holds because $Q(A/P)$ is flat as a left $A/N$-module (by 
\cite{Pro, Ch. V, Th. 2.5} or \cite{Row, Th. 3.2.27} the direct product of the 
rings $Q(A/P')$ with $P'$ running over the set of minimal prime ideals of $A$ 
is a classical quotient ring of the finitely generated semiprime PI algebra 
$A/N$). In general the sum in the definition of $\rho_P(M)$ does not change 
when a chosen Loewy series is replaced with any its refinement. Since any two 
Loewy series have refinements with isomorphic factors, it is immediate that 
the sum does not depend on the choice.

The reduced rank $\rho(M)$ of $M$ is the sum $\sum\rho_P(M)$ over all minimal 
primes $P$. However, we will need the values $\rho_P(M)$ only at primes 
$P\in\Om(A)$, while there could exist minimal primes outside this set.

\proclaim
Lemma 1.3.
Any $A$-module $M$ has a finite chain of submodules 
$M=M_0\sps M_1\sps\cdots\sps M_n=0$ in which each factor is annihilated by 
a prime ideal of $A$.
\endproclaim

\Proof.
This is a consequence of the fact that $P_1\cdots P_n=0$ for a suitable sequence 
$P_1,\ldots,P_n\in\Spec A$, as we have noted already in the proof of Lemma 1.2.
\endproof

\proclaim
Lemma 1.4.
For a right $A$-module $M$ we have $M\in\T_r$ if and only if $\rho_P(M)=0$ for 
each $P\in\Om(A)$.
\endproclaim

\Proof.
Take a chain of submodules $M=M_0\sps M_1\sps\cdots\sps M_n=0$ in which each 
factor $M_{i-1}/M_i$ is annihilated by a prime ideal, say $P_i$, of $A$. Such 
a chain exists by Lemma 1.3, and it is a Loewy series since $N\sbs P_i$ for 
each $i=1,\ldots,n$.

For $P\in\Om(A)$ the equality $\rho_P(M)=0$ holds if and only if
$$
M_{i-1}/M_i\ot_AQ(A/P)=0
$$
for each $i$. If $P_i\ne P$ for some $i$, then the latter equality is always 
true since the image of $P_i$ in the Goldie ring $A/P$ is a nonzero ideal, and 
therefore there exists $s\in P_i$ such that the coset $s+P$ is a regular 
element of $A/P$, while at the same time $s$ annihilates $M_{i-1}/M_i$. In 
other words, only those indices $i$ for which $P_i=P$ have to be checked.

Thus $\rho_P(M)=0$ for each $P\in\Om(A)$ if and only if 
$M_{i-1}/M_i\ot_AQ(A/P_i)=0$ for each $i$ with $P_i\in\Om(A)$. By Lemma 1.1 
this is equivalent to the condition that $M_{i-1}/M_i\in\T_r$ for all such 
indices $i$. We also know that $M_{i-1}/M_i\in\T_r$ for all $i$ with 
$P_i\notin\Om(A)$, and the desired conclusion follows from the fact that 
$M\in\T_r$ if and only if $M_{i-1}/M_i\in\T_r$ for each $i=1,\ldots,n$.
\endproof

\setitemsize(ii)
\proclaim
Lemma 1.5.
Let $M'$, $M''$ be two submodules of a right $A$-module $M$. Suppose that 
$\rho_P(M)<\infty$ for each $P\in\Om(A)$. Then\/{\rm:}

\item(i)
\ $M/M'\in\T_r$ if and only if $\rho_P(M)=\rho_P(M')$ for each $P\in\Om(A)$.  

\item(ii)
\ If $M/M'\in\T_r$ and $M''\cong M',$ then $M/M''\in\T_r$ too.

\endproclaim

\Proof.
Let us build a Loewy series for $M$ by first taking the preimages in $M$ of 
the terms of any chosen Loewy series for $M/M'$, and then adding the remaining 
terms from a Loewy series for $M'$. The factors $M_{i-1}/M_i$ of the Loewy 
series thus obtained will be those of the two Loewy series used in the 
construction. From this it is clear that
$$
\rho_P(M)=\rho_P(M')+\rho_P(M/M')
$$
for each minimal prime $P$ of $A$. Finiteness of $\rho_P(M)$ implies that 
$\rho_P(M')<\infty$ too and $\rho_P(M/M')=\rho_P(M)-\rho_P(M')$. Lemma 1.4 
tells us that $M/M'\in\T_r$ if and only if $\rho_P(M/M')=0$ for each 
$P\in\Om(A)$, which amounts to (i).

Similarly, $M/M''\in\T_r$ if and only if $\rho_P(M)=\rho_P(M'')$ for each 
$P\in\Om(A)$. If $M''\cong M'$, then $\rho_P(M'')=\rho_P(M')$ for all minimal 
primes $P$ of $A$, and (ii) is immediate.
\endproof

\proclaim
Corollary 1.6.
Let $\al\in\End_AM$. If\/ $\Ker\al=0$ and $\rho_P(M)<\infty$ for each 
$P\in\Om(A)$ then $M/\al(M)\in\T_r$.
\endproclaim

\Proof.
Take $M'=M$ and $M''=\al(M)\cong M$ in Lemma 1.5.
\endproof

\proclaim
Lemma 1.7.
Let $M$ be a right $A$-module{\rm,} $M'$ its submodule such that $M/M'\in\T_r$. 
Then $\rho_P(M)=\rho_P(M')$ for each $P\in\Om(A)$. In particular{\rm,} 
$\rho_P(M)<\infty$ if and only if $\rho_P(M')<\infty$.
\endproclaim

\Proof.
This is immediate from the equality $\rho_P(M)=\rho_P(M')+\rho_P(M/M')$ which 
we have observed in the proof of Lemma 1.5. Indeed, $\rho_P(M/M')=0$ for each 
$P\in\Om(A)$ by Lemma 1.4.
\endproof

Lemma 1.9 below serves as a means to verify that $\rho_P(M)<\infty$. Its proof 
makes use of Lemma 1.8. A bimodule $M$ is said to be \emph{$k$-central} if the 
two module structures on $M$ are extensions of the same $k$-vector space 
structure, i.e., $av=va$ for all $a\in k$ and $v\in M$.

\proclaim
Lemma 1.8.
Let $M$ be a $k$-central $(F,E)$-bimodule where $F$ and $E$ are extension 
fields of the base field $k$. Suppose that $F$ is a finitely generated 
extension of $k$ with
$$
\trdeg F/k=\trdeg E/k.
$$
If\/ $\dim_FM<\infty,$ then\/ $\dim_EM<\infty$.
\endproclaim

\Proof.
We proceed by induction on $\dim_FM$. If $\dim_FM=0$, then $M=0$. Suppose that 
$\dim_FM>0$. Consider the algebraic closure $K$ of $F$. The action of $E$ on 
$M$ gives an embedding $\,E\sbs\End_FM\sbs\End_K(K\ot_FM)$. In this way $E$ 
becomes an algebra of linear transformations of the finite-dimensional vector 
space $K\ot_FM$ over the field $K$. Since $E$ is commutative, there is a common 
eigenvector $v\in K\ot_FM$ for all elements of $E$, so that $vE\sbs Kv$. 

We have $v\in L\ot_FM$ for some finite field extension $L$ of $F$ contained in 
$K$. Now $L\ot_FM$ is an $(L,E)$-bimodule in which $M'=Lv$ is a subbimodule since
$$
Lv=Kv\cap(L\ot_FM)\sbs K\ot_FM.
$$
There is an embedding of fields $\la:E\to L$ such that $\la(a)v=va$ for $a\in E$. 
Since $M$ is $k$-central, we have $\la(a)=a$ for all $a\in k$. Note that $L$ 
is a finitely generated extension of $k$ with $\,\trdeg L/k=\trdeg F/k=\trdeg E/k$. 
It follows that $\,\trdeg L/\la(E)=0$. Thus $L$ is a finitely generated 
algebraic extension of its subfield $\la(E)$. This extension has to be finite. 
Hence
$$
\dim_EM'=[L:\la(E)]<\infty.
$$
Consider the $(L,E)$-bimodule $M''=(L\ot_FM)/M'$. Since 
$$
\dim_LM''<\dim_L(L\ot_FM)=\dim_FM,
$$
we may assume by the induction hypothesis that $\,\dim_EM''<\infty\mskip1mu$. Then
$$
\dim_EM\le\dim_E(L\ot_FM)=\dim_EM'+\dim_EM''<\infty,
$$
and we are done.
\endproof

\proclaim
Lemma 1.9.
Let $M$ be a $k$-central $(F,A)$-bimodule where $F$ is a finitely generated 
extension field of the base field $k$. Suppose that\/ $\trdeg F/k=d(A)$. 
If\/ $\dim_FM<\infty,$ then $\rho_P(M)<\infty$ for each $P\in\Om(A)$. 
\endproclaim

\Proof.
There is a chain of $A$-submodules $M=M_0\sps M_1\sps\cdots\sps M_n=0$ in 
which each factor is annihilated by a prime ideal of $A$. In the construction 
of such a chain we may take $M_i=M_{i-1}P_i$ with $P_i\in\Spec A$ for each 
$i=1,\ldots,n$. Then each $M_i$ is a subbimodule of $M$. Since 
$\rho_P(M)=\sum\rho_P(M_{i-1}/M_i)$, it suffices to consider the case 
$n=1$ when $M$ is annihilated by $P_1$.

Let $P\in\Om(A)$. If $P\ne P_1$, then $M\ot_AQ(A/P)=0$, and $\rho_P(M)=0$. 
Suppose that $P=P_1$. Then $M$ is a right $A/P$-module. Its torsion submodule 
$M'$ consists of all elements whose annihilator in $A/P$ contains a regular 
element of that ring. We have $M'\ot_AQ(A/P)=0$. Clearly $M'$ is a subbimodule 
of $M$.

The factor bimodule $M''=M/M'$ is torsionfree as a right $A/P$-module. Each 
regular element of $A/P$ acts on $M''$ as an injective $F$-linear 
transformation. Since $\,\dim_FM''<\infty$, this transformation is invertible. 
Hence $M''$ is a right $Q(A/P)$-module in a natural way. The center $E$ of 
$Q(A/P)$ is an extension field of $k$. By \cite{Kr-L, Th. 10.5} or 
\cite{Row, Prop. 6.3.40}
$$
\trdeg E/k=\GKdim A/P=d(A)=\trdeg F/k.
$$
Lemma 1.8 applied to the $(F,E)$-bimodule $M''$ yields $\,\dim_EM''<\infty$. 
Hence $M''$ is a finitely generated $Q(A/P)$-module, and so is 
$\,M\ot_AQ(A/P)\cong M''$.
\endproof

Hereditary torsion theories on the categories of modules are characterized by 
means of certain filters of one-sided ideals called Gabriel topologies 
\cite{St, Ch. VI, Th. 5.1}. The filter $\G_r$ corresponding to the torsion 
class $\T_r$ consists of all right ideals $I$ of the algebra $A$ such that 
$A/I\in\T_r$. It is a set of right ideals satisfying the following properties 
of a \emph{right Gabriel topology}:

\setitemsize(T4)
\item(T1)
if $J\in\G_r$ and $I$ is a right ideal such that $J\sbs I$ then $I\in\G_r$,

\item(T2)
if $I,J\in\G_r$ then $I\cap J\in\G_r$,

\item(T3)
if $I\in\G_r$ then $(I:a)\in\G_r$ for each $a\in R$,

\item(T4)
if $J\in\G_r$ and $I$ is a right ideal such that $(I:a)\in\G_r$ for all 
$a\in J$ then $I\in\G_r$,

\noindent
where $(I:a)=\{x\in A\mid ax\in I\}$.

A right $A$-module $M$ belongs to the torsion class $\T_r$ if and only if each 
element of $M$ is annihilated by a right ideal in $\G_r$. We will use the term 
\emph{$\G_r$-torsion module} for a module in this torsion class.

An arbitrary right $A$-module $M$ has a largest $\G_r$-torsion submodule. We say 
that $M$ is \emph{$\G_r$-torsionfree} if it contains no nonzero $\G_r$-torsion 
submodules. Such modules form the torsionfree class of the torsion theory 
considered here. By the definition of $\T_r$ any nonzero $\G_r$-torsion module 
contains a nonzero submodule of Gelfand-Kirillov dimension less than $d(A)$. 
Therefore a right $A$-module $M$ is $\G_r$-torsionfree if and only if 
$\,\GKdim M'\ge d(A)\,$ for each nonzero submodule $M'\sbs M$.

The previous definitions were concerned with right $A$-modules and right 
ideals of $A$. Denote by $\T_l$ the torsion class of left $A$-modules defined 
similarly to $\T_r$, and let $\G_l$ be the left Gabriel topology determined by 
the left ideals $I$ of $A$ such that $A/I\in\T_l$. We will speak in a similar 
way about \emph{$\G_l$-torsion} and \emph{$\G_l$-torsionfree} left $A$-modules. 
The reduced ranks $\rho_P(M)$ make sense also for left $A$-modules.

\setitemsize(a)
\proclaim
Lemma 1.10.
For a two-sided ideal $I$ of $A$ the following conditions are equivalent\/{\rm:}

\item(a)
$\ I\in\G_l,$

\item(b)
$\ I\in\G_r,$

\item(c)
$\ \GKdim A/P<d(A)$ for each prime ideal $P$ of $A$ such that $I\sbs P$.

\endproclaim

\Proof.
By Lemma 1.3 the factor algebra $A/I$ has a finite chain of right ideals in 
which each factor is annihilated by a prime ideal of $A/I$, i.e., each factor 
is annihilated by a prime ideal of $A$ containing $I$. Property (c) implies 
that each of these factors has Gelfand-Kirillov dimension less than $d(A)$, 
and therefore $A/I\in\T_r$, i.e., $I\in\G_r$. This shows that (c)$\Rar$(b).

Conversely, if $A/I\in\T_r$, then $A/P\in\T_r$ for each prime ideal $P$ of $A$ 
containing $I$, but this is equivalent to (c) in view of Lemma 1.1. Thus 
(b)$\,\Lrar\,$(c). By symmetry we also have (a)$\,\Lrar\,$(c).
\endproof

\section
2. The symmetric quotient ring

Each Gabriel topology on a ring $R$ gives rise to a localization of $R$ 
\cite{St, Ch. IX}. We are going to apply this construction to the Gabriel 
topology $\G_r$ introduced in section 1. We continue to assume that $A$ is a 
finitely generated PI algebra. Furthermore, we will need several additional 
assumptions:

\setitemsize(A3)
\item(A1)
\ $A$ satisfies the ACC on right and left annihilators,

\item(A2)
\ $\GKdim I\ge d(A)$ for each nonzero one-sided ideal $I$ of $A$,

\item(A3)
\ $\rho_P(A_A)<\infty$ and $\rho_P({}_AA)<\infty$ for each $P\in\Om(A)$.

\smallskip
The subscripts in $A_A$ and ${}_AA$ indicate that $A$ is regarded as a module 
over itself with respect to either right or left multiplications, respectively. 
Condition (A2) should be understood as the bound on the Gelfand-Kirillov 
dimension of all right ideals regarded as right $A$-modules and all left 
ideals regarded as left $A$-modules. It means precisely that $A_A$ is 
$\G_r$-torsionfree and ${}_AA$ is $\G_l$-torsionfree.

The \emph{right quotient ring} of $A$ with respect to $\G_r$ is defined as the 
direct limit
$$
Q_r(A)=\limdir_{I\in\G_r}\Hom_A(I,A).
$$
Each element of $Q_r(A)$ is represented by a right $A$-linear map $\al:I\to A$ 
with a right ideal $I\in\G_r$, and another right $A$-linear map $\be:J\to A$ 
with $J\in\G_r$ represents the same element of $Q_r(A)$ if and only if $\be$ 
agrees with $\al$ on a smaller right ideal from the filter $\G_r$. The 
multiplication in $Q_r(A)$ is induced by the composition of maps. The 
algebra $A$ is identified with the subring of $Q_r(A)$ consisting of all 
elements represented by left multiplications in $A$.

For each $q\in Q_r(A)$ there exists $I\in\G_r$ such that $qI\sbs A$, and 
$qI\ne0$ whenever $q\ne0$. Indeed, if $q$ is represented by $\al:I\to A$ then 
$qx=\al(x)$ for all $x\in I$. It follows that $Q_r(A)$ is $\G_r$-torsionfree 
as a right $A$-module, and each nonzero right $A$-submodule of $Q_r(A)$ has 
nonzero intersection with $A$. (This implies that $Q_r(A)$ is a rational 
extension of $A$ in the category of right $A$-modules, and so $Q_r(A)$ is a 
subring of the maximal right quotient ring of $A$, as discussed in 
\cite{Goo}.)

Now we introduce the \emph{symmetric quotient ring} of $A$. Put
$$
Q(A)=\{q\in Q_r(A)\mid\hbox{ there exists }\,J\in\G_l\,
\hbox{ such that }\,Jq\sbs A\}.
$$
We will write $Q_r=Q_r(A)$ and $Q=Q(A)$ for short, and the notation $Q_A$ will 
be used when $Q$ is regarded as a right $A$-module.

\proclaim
Lemma 2.1.
The set $Q$ is a subring of $Q_r$ containing $A$. It is $\G_l$-torsionfree as 
a left $A$-module and $\G_r$-torsionfree as a right $A$-module.

For each $q\in Q$ there exist one-sided ideals $J\in\G_l$ and $I\in\G_r$ such 
that $Jq\sbs A$ and $qI\sbs A$. In other words{\rm,} the $A$-bimodule $Q/A$ 
is $\G_l$-torsion and $\G_r$-torsion.
\endproclaim

\Proof.
By the definition above $Q$ consists of all $q\in Q_r$ such that the coset 
$q+A$ is a $\G_l$-torsion element of the left $A$-module $Q_r/A$. By the 
axioms of Gabriel topologies the $\G_l$-torsion elements form a submodule. 
This means that $Q$ is a left $A$-submodule of $Q_r$ such that $Q/A$ is 
$\G_l$-torsion and $Q_r/Q$ is $\G_l$-torsionfree.

Given $x\in Q$, we thus have $Ax\sbs Q$. It follows that $Qx/(Qx\cap Q)$ is an 
epimorphic image of the left $A$-module $Q/A$. Hence $Qx/(Qx\cap Q)$ is 
$\G_l$-torsion. On the other hand, this module is $\G_l$-torsionfree since it 
embeds in $Q_r/Q$. This is possible only when $Qx/(Qx\cap Q)=0$, i.e., $Qx\sbs Q$. 
Hence $Q$ is closed under the multiplication in $Q_r$, and so $Q$ is indeed a 
subring.

Suppose that $Jq=0$ for some $q\in Q_r$ and $J\in\G_l$. There exists 
$I\in\G_r$ such that $qI\sbs A$. Since $JqI=0$ and $A$ is $\G_l$-torsionfree, 
we deduce that $qI=0$, but then $q=0$. This shows that $Q_r$ is not only 
$\G_r$-torsionfree, but also $\G_l$-torsionfree. Hence so too is $Q$.

The last assertion in the lemma is clear from the definitions of $Q_r$ and $Q$.
\endproof

\Remark 1.
Consider the left quotient ring $Q_l$ of $A$ defined with respect to the left 
Gabriel topology $\G_l$. Then $Q$ is canonically isomorphic with the subring 
of $Q_l$ consisting of all elements $q\in Q_l$ such that $qI\sbs A$ for some 
$I\in\G_r$. In this sense the construction of $Q$ is left-right symmetric.
\endremark

The initial reason which forces us to work with the subring $Q$, rather than 
with the whole ring $Q_r$, is that for $Q_r$ we do not have the analog of the 
next lemma. The ACC on annihilators will be important in proving that $Q$ is 
semiprimary.

For a ring $R$ and its subset $X$ we denote by $\,\rann_RX\,$ and 
$\,\lann_RX\,$ the right and left annihilators of $X$ in $R$.

\proclaim
Lemma 2.2.
The ring $Q$ satisfies the ACC on right and left annihilators.
\endproclaim

\Proof.
Given $q\in Q$, there exists $J\in\G_l$ such that $Jq\sbs A$. Since $Q$ is 
$\G_l$-torsionfree, we have $\,\rann_Qq=\rann_QJq\mskip1mu$. It follows that 
for each subset $X\sbs Q$ there exists $X'\sbs A$ such that $\,\rann_QX=\rann_QX'$.

Let $p\in Q$, and let $I\in\G_r$ be such that $p\mskip1mu I\sbs A$. Since $Q$ 
is $\G_r$-torsionfree, the equality $Xp=0$ is equivalent to $Xp\mskip1mu I=0$. 
It follows that
$$
p\in\rann_QX\quad\hbox{if and only if}\quad p\mskip1mu I\sbs A\cap\rann_QX.
$$

If $K_1\sbs K_2\sbs\cdots$ is an ascending chain of annihilator right ideals 
of $Q$, then their intersections with $A$ form an ascending chain of 
annihilator right ideals of $A$. Since $A$ satisfies the ACC on right 
annihilators, there is an integer $n>0$ such that $A\cap K_i=A\cap K_n$ for all 
$i>n$. But each annihilator right ideal of $Q$ is completely characterized by 
its intersection with $A$, as explained in the preceding paragraph. Therefore 
$K_i=K_n$ for all $i>n$. This shows that $Q$ satisfies the ACC on right 
annihilators, and the case of left annihilators is similar.
\endproof

\proclaim
Lemma 2.3.
Let $\,I=uQ+\rann_Qu\,$ where $u\in Q$ is an element satisfying 
$$
\rann_Qu=\rann_Qu^2.
$$
Then $I$ is a right ideal of $Q$ such that $Q/I$ is $\G_r$-torsion as a right 
$A$-module. We have $\,A\cap I\in\G_r\,$ and{\rm,} moreover{\rm,} 
$\,A\cap uQ+A\cap\,\rann_Qu\in\G_r\,$.
\endproclaim

\Proof.
Put $M=Q/\rann_Qu$. If $x\in Q$ is any element such that $ux\in\rann_Qu$, then 
$u^2x=0$, whence $x\in\rann_Qu$ by the assumption about $u$. This shows that 
the left multiplication by $u$ induces an injective endomorphism $\al$ of $M$ 
as a right $Q$-module, and therefore also as a right $A$-module. 

As we know, $Q/A$ is $\G_r$-torsion, i.e., $Q/A\in\T_r$. By Lemma 1.7 
assumption (A3) implies that $\rho_P(Q_A)<\infty$ for each $P\in\Om(A)$. Since 
$\rho_P(M)\le\rho_P(Q_A)$, we have $\rho_P(M)<\infty$ for each $P\in\Om(A)$. 
Now the first conclusion follows from Corollary 1.6 since $\,Q/I\cong M/\al(M)$.

The right $A$-module $A/(A\cap I)$ is $\G_r$-torsion since it embeds in $Q/I$. 
This means that $A\cap I\in\G_r$. Next, put
$$
I_1=A\cap uQ\quad{\rm and}\quad I_2=A\cap\mskip1mu\rann_Qu.
$$
For each right ideal $K$ of $Q$ the right $A$-module $K/(A\cap K)$ is 
$\G_r$-torsion since it embeds in $Q/A$. In particular, both $uQ/I_1$ and 
$(\rann_Qu)/I_2$ are $\G_r$-torsion. Hence so is $I/(I_1+I_2)$, being the sum 
of epimorphic images of the previous two $A$-modules. Since $Q/I$ is also 
$\G_r$-torsion, so too is $Q/(I_1+I_2)$. It follows that $A/(I_1+I_2)\in\T_r$, 
which yields $\,I_1+I_2\in\G_r$.
\endproof

\proclaim
Lemma 2.4.
If $I$ is a right ideal of $Q$ such that $A\cap I\in\G_r,$ then each right 
$A$-linear map $f:I\to Q$ is induced by a left multiplication in the larger 
ring $Q_r$.
\endproclaim

\Proof.
Put $I_0=A\cap I$ and $I'=A\cap f^{-1}(A)$. Both are right ideals of $A$ and 
$I'\sbs I_0$. The right $A$-module $I_0/I'$ is $\G_r$-torsion since it embeds 
in $Q/A$ by means of the map induced by $f$. Since $A/I_0$ is also 
$\G_r$-torsion, so is $A/I'$, i.e., $I'\in\G_r$.

The right $A$-linear map $f|_{I'}:I'\to A$ represents an element $q\in Q_r$ 
such that $qx=f(x)$ for all $x\in I'$. It remains to show that $qx=f(x)$ for 
all $x\in I$. 

The right $A$-linear map $h:I\to Q_r$ defined by the rule $h(x)=qx-f(x)$ for 
$x\in I$ vanishes on $I'$. Hence $h$ factors through $I/I'$. Note that the 
right $A$-module $I/I'$ is a submodule of $Q/I'$, and the latter is an 
extension of $A/I'$ by $Q/A$. It follows that $I/I'$ is $\G_r$-torsion. 
Since $Q_r$ is $\G_r$-torsionfree, we have $\Hom_A(I/I',Q_r)=0$. Therefore 
$h=0$, and we are done.
\endproof

\proclaim
Lemma 2.5.
Given $u\in Q$ satisfying $\,\rann_Qu=\rann_Qu^2$ and $\,\lann_Qu=\lann_Qu^2,$ 
there exist two elements $e,v\in Q$ such that
$$
eu=ue=u,\quad ev=ve=v,\quad uv=vu=e,\quad e^2=e.
$$
\endproclaim

\Proof.
In the ring $Q$ consider its right and left ideals
$$
I=uQ+\rann_Qu,\qquad J=Qu+\lann_Qu.
$$
By Lemma 2.3 $\,A\cap I\in\G_r$.

The assumption about $u$ implies that the sum $uQ+\rann_Qu$ is direct. Hence 
there is a right $Q$-linear map $f:I\to Q$ such that $f(x)=x$ for all $x\in uQ$ 
and $f(x)=0$ for all $x\in\rann_Qu$. By Lemma 2.4 $f$ is the restriction to $I$ 
of the left multiplication by some element $e\in Q_r$. Then
$$
eu=u\quad{\rm and}\quad\rann_Qu\sbs\rann_Qe.
$$
Since $(ue-u)I=0$ and $(e^2-e)I=0$, it follows that $ue=u$ and $e^2=e$, again 
by $\G_r$-torsionfreeness of $Q_r$. Indeed, $\,\lann_{Q_r}I=0\,$ since 
$\,A\cap I\in\G_r$.

We have to show that $e\in Q$. For this we have to find a left ideal $J'$ of 
$A$ with the properties that $J'\in\G_l$ and $J'e\sbs A$. Let us take
$$
J'=J_1+J_2\quad\hbox{where $\,J_1=A\cap Qu\,$ and $\,J_2=A\cap\lann_Qu$}.
$$
By the left hand version of Lemma 2.3 applied to $J$ we have $J'\in\G_l$.

If $y\in\lann_Qu$, then $yeu=yu=0$. Since $ex=0$ for all $x\in\rann_Qu$, we 
deduce that $yeI=0$, whence $ye=0$ by $\G_r$-torsionfreeness of $Q_r$. In 
particular, $ye=0$ for all $y\in J_2$. For $t\in Qu$ we have $te=t$ since 
$ue=u$. In particular, $te=t\in A$ for all $t\in J_1$. We see that $J'e\sbs A$, 
and the inclusion $e\in Q$ is thus proved.

Since $\rann_Qu\sbs\rann_Qe$, there is a right $Q$-linear map $g:I\to Q$ such 
that $g(u)=e$ and $g(x)=0$ for all $x\in\rann_Qu$. By Lemma 2.4 $g$ is the 
restriction to $I$ of the left multiplication by some element $v\in Q_r$. Then
$$
vu=e\quad{\rm and}\quad\rann_Qu\sbs\rann_Qv.
$$
Observing that $\,(uv-e)I=0$, $\,(ev-v)I=0$, $\,(ve-v)I=0\,$, we get the 
equalities $\,uv=e\,$ and $\,ev=ve=v$.

Consider the left ideals $J_1,J_2,J'$ of $A$ defined earlier, and put 
$J''=K_1+J_2$ where $K_1=A\cap J'u\sbs J_1$. The left $A$-module $J_1/K_1$ 
embeds in $Qu/J'u$, which is an epimorphic image of $Q/J'\in\T_l$. Hence all 
these modules are $\G_l$-torsion, and so too is $J'/J''$. Since $J'\in\G_l$, 
we get $J''\in\G_l$.

For all $y\in J_2$ we have $yvI=0$ since $yvu=ye=0$ and $vx=0$ for all 
$x\in\rann_Qu$. Hence $J_2v=0$ by $\G_r$-torsionfreeness of $Q_r$. Also 
$K_1v\sbs A$ since $J'uv=J'e\sbs A$. It follows that $J''v\sbs A$, and 
therefore $v\in Q$.
\endproof

A ring is called \emph{semiprimary} if its Jacobson radical is nilpotent and 
the factor ring by the Jacobson radical is semisimple Artinian.

\proclaim
Lemma 2.6.
The ring $Q$ is semiprimary.
\endproclaim

\Proof.
Since the ring $Q$ satisfies the ACC on right and left annihilators by Lemma 
2.2, all its nil subrings are nilpotent \cite{Her-S64, Th. 1}. Hence $Q$ has a 
nilpotent ideal $N$ which contains every nil right ideal of $Q$.

We claim that each right ideal $I$ of $Q$ has the form $I=eQ+K$ where $e\in Q$ 
is an idempotent and $K$ a nil right ideal of $Q$, so that $K\sbs N$. In order 
to prove this consider the set
$$
X=\{u\in I\mid\rann_Qu=\rann_Qu^2\ \ {\rm and}\ \ \lann_Qu=\lann_Qu^2\}.
$$
The lattices of annihilator right and left ideals of $Q$ are antiisomorphic to 
each other. Therefore the ACC on left annihilators implies that $Q$ also 
satisfies the DCC on right annihilators. In particular, the set 
$\{\rann_Qx\mid x\in X\}$ has a minimal element. Pick $u\in X$ such that 
$\,\rann_Qu\,$ is a minimal element of that set of right ideals.

By Lemma 2.5 there is an idempotent $e\in Q$ such that $uQ=eQ$ and $Qu=Qe$. 
The first of these two equalities shows that $e\in I$ since $u\in I$. The 
second equality implies that $\,\rann_Qu=\rann_Qe=(1-e)Q\,$. Since 
$Q=eQ\oplus(1-e)Q$, we get
$$
I=eQ\oplus K\quad\hbox{where $K=I\cap(1-e)Q$}.
$$
We will prove that $K$ is nil. Suppose on the contrary that $K$ contains a 
nonnilpotent element $y$. The right ideals $\rann_Qy^i$, $\,i=1,2,\ldots\,$, 
form an ascending chain which has to stabilize by the ACC. Replacing $y$ with 
$y^n$ for sufficiently large $n$, we may assume therefore that 
$\rann_Qy^2=\rann_Qy$.

Put $t=e+y\in I$. Since $y\in K$, we have $eQ\cap yQ=0$. Hence 
$$
\rann_Qt=\rann_Qe\cap\rann_Qy\sbs(1-e)Q\,.
$$
Since $ey=0$, but $ty=y^2\ne0$, the last inclusion is proper. Next, 
$t^2=e+ye+y^2$. If $t^2q=0$ for some $q\in Q$, then $eq=0$ since $eQ\cap 
yQ=0$, and we must also have $y^2q=0$. It follows that
$$
\rann_Qt^2=\rann_Qe\cap\rann_Qy^2=\rann_Qt,
$$
and therefore $\,\rann_Qt^i=\rann_Qt\,$ for all $i>0$. Since the ascending chain 
of left ideals $\,\lann_Qt^i$, $\,i=1,2,\ldots\,$, stabilizes, we have $t^m\in X$ 
for sufficiently large $m$. On the other hand, $\,\rann_Qt^m\,$ is properly 
contained in $\,\rann_Qu=(1-e)Q$, which contradicts the choice of $u$.

Thus our claim about the right ideals of $Q$ has been proved. It follows that 
each right ideal of the factor ring $Q/N$ is generated by an idempotent. Hence 
$Q/N$ is semisimple Artinian, and $N$ is the Jacobson radical of $Q$.  
\endproof

\proclaim
Lemma 2.7.
If $A'$ is any finitely generated PI algebra and $M$ a nonzero $A'$-module{\rm,} 
then $M$ has a nonzero element annihilated by a prime ideal of $A'$.
\endproclaim

\Proof.
This is immediately clear from Lemma 1.3.
\endproof

The next statement repeats \cite{Sk99, Lemma 3.8}:

\setitemsize(a)
\proclaim
Lemma 2.8.
Let $R$ be a semiprime right Goldie subring of a semisimple Artinian ring $S$. 
Suppose that $\G$ is a set of right ideals of $R$ with the two properties\/{\rm:}

\item(a)
\ $\lann_SI=0$ for each $I\in\G,$

\item(b)
\ for each $x\in S$ there exists $I\in\G$ such that $xI\sbs R$.

Then $S$ is a classical right quotient ring of $R,$ and each right ideal 
$I\in\G$ contains a regular element of $R$.
\endproclaim

When we apply notions of the $\G_l$-torsion and $\G_r$-torsion theories to 
$Q$-modules we always regard $Q$-modules as $A$-modules.

\setitemsize(A4)
\proclaim
Proposition 2.9.
In addition to {\rm(A1), (A2), (A3)} assume another condition\/{\rm:}

\item(A4)
each right $Q$-module is $\G_r$-torsionfree and each left $Q$-module is 
$\G_l$-torsionfree.

Then $Q$ is a right and left Artinian classical quotient ring of $A$.
\endproclaim

\Proof.
In several arguments below we follow the proof of \cite{Sk99, Prop. 3.9} with 
some modifications. By Lemma 2.6 $Q$ is semiprimary. Hence $Q$ has finitely 
many maximal ideals, and all factor rings by those ideals are simple Artinian.

\proclaim
Claim 1.
If $M$ is any maximal ideal of $Q,$ then $M\cap A$ is a prime ideal of $A$.
\endproclaim

There is an embedding of factor rings $A/(M\cap A)\hrar Q/M$. Viewing $Q/M$ as 
a right $A/(M\cap A)$-module, we deduce from Lemma 2.7 that $Q/M$ has a 
nonzero element whose right annihilator contains a prime ideal of 
$A/(M\cap A)$. In other words, there exist $y\in Q$ and $P\in\Spec A$ such 
that $y\notin M$, $M\cap A\sbs P$ and $yP\sbs M$.

For each $I\in\G_r$ put
$$
L(I)=\{x\in Q\mid xIP\sbs M\}.
$$
It is clear that $L(I)$ is a left ideal of $Q$ and $M\sbs L(I)$, so that 
$L(I)/M$ is a left ideal of the factor ring $Q/M$. Note that for 
$I_1,I_2\in\G_r$ there exists $I_3\in\G_r$ such that $L(I_3)$ contains both 
$L(I_1)$ and $L(I_2)$. Indeed, $I_3=I_1\cap I_2$ will do. Since the simple 
ring $Q/M$ is Artinian, it follows that the set of left ideals 
$\{L(I)\mid I\in\G_r\}$ has a largest element $L_0$. Thus $L(I)\sbs L_0$ 
for each $I\in\G_r$.

Note that $L(A)=\{x\in Q\mid xP\sbs M\}$ is a right $A$-submodule of $Q$. Let 
$x\in L(A)$ and $q\in Q$. By Lemma 2.1 there exists $I\in\G_r$ such that 
$qI\sbs A$. Then
$$
xqIP\sbs xAP\sbs xP\sbs M,
$$
yielding $xq\in L(I)\sbs L_0$. This shows that $L(A)Q\sbs L_0$. Note that 
$L(A)Q$ is a two-sided ideal of $Q$ containing $M$. But $L(A)\not\sbs M$ by 
the choice of $P$, whence $L(A)Q\not\sbs M$. Since $Q/M$ is a simple ring, we 
must have $L(A)Q=Q$. It follows that $L_0=Q$ too.

Now pick $I\in\G_r$ such that $L_0=L(I)$. Since $1\in L_0$, we get $IP\sbs M$. 
Then also $I'P\sbs M$ where $I'=AI$ is a two-sided ideal of $A$ containing $I$. 
Since $I'\in\G_r$, we have $I'\in\G_l$ by Lemma 1.10. This shows that the image 
of $P$ in $Q/M$ is a $\G_l$-torsion left $A$-submodule of $Q/M$. Since $Q/M$ 
is $\G_l$-torsionfree by condition (A4), we get $P\sbs M$. Hence $P=M\cap A$, 
and Claim 1 is thus proved.

\proclaim
Claim 2.
Denote by $N$ the prime radical of $A$ and by $J$ the Jacobson radical of $Q$. 
Then $N=J\cap A$ and $J=NQ$. 
\endproclaim

Let $M_1,\ldots,M_k$ be all the maximal ideals of $Q$. Then $J=\bigcap M_i$. 
By Claim 1 $P_i=M_i\cap A$ is a prime ideal of $A$ for each $i$. We have 
$J\cap A=N'$ where $N'=\bigcap P_i$. Since $J$ is nilpotent, $N'$ is a 
nilpotent ideal of $A$. Hence $N'$ is contained in each prime ideal of $A$, 
and therefore $N'=N$.

The equality $N=J\cap A$ implies that $N\sbs J$ and $J/N$ is a $\G_r$-torsion 
right $A$-module. Then $NQ\sbs J$ and $J/NQ$ is $\G_r$-torsion as well. Since 
$J/NQ$ is a right $Q$-module, it is $\G_r$-torsionfree by (A4). It follows 
that $J/NQ=0$, i.e., $J=NQ$.

\proclaim
Claim 3.
The factor ring $Q/J$ is a classical quotient ring of $A/N$.
\endproclaim

Put $S=Q/J$ and $R=A/N$. The ring $S$ is semisimple Artinian, while $R$ is a 
finitely generated semiprime PI algebra. By \cite{Pro, Ch. V, Th. 2.5} $R$ is 
semiprime right and left Goldie. Let $\pi:Q\to S$ be the canonical 
homomorphism. By Claim 2 $R$ is identified with the subring $\pi(A)$ of $S$.

Consider the set $\{\pi(I)\mid I\in\G_r\}$ of right ideals of $R$. By (A4) $S$ 
is $\G_r$-torsionfree as a right $A$-module, which means that $\lann_S\pi(I)=0$ 
for each $I\in\G_r$. If $q\in Q$, then $qI\sbs A$, and therefore 
$\pi(q)\pi(I)\sbs\pi(A)$, for some $I\in\G_r$. Now Claim 3 follows from Lemma 2.8.

\medbreak
Denote by $\C$ the set of all elements $u\in A$ which are regular modulo $N$, 
i.e., whose images $\pi(u)$ in the ring $\pi(A)\cong A/N$ are regular elements 
of that ring. If $u\in\C$, then $\pi(u)$ is invertible in $Q/J$ by Claim 3; 
this implies that $u$ is invertible in $Q$ since $J$ is the Jacobson radical 
of $Q$. Thus all elements of $\C$ are regular in $A$.

Conversely, suppose that $u$ is any regular element of $A$. Then $\rann_Qu$ is 
a right ideal of $Q$ which has zero intersection with $A$. This entails 
$\rann_Qu=0$. Similarly, $\lann_Qu=0$, and so $u$ is regular in $Q$. Such an 
element satisfies the assumptions of Lemma 2.5. Hence $uQ=eQ$ and $Qu=Qe$ for 
some idempotent $e\in Q$. But then $(1-e)u=0$, which forces $e=1$. We conclude 
that $u$ is invertible in $Q$ (more generally, it is known that any regular 
element of a semiprimary ring is invertible). Since $\pi(u)$ is then 
invertible in $Q/J$, we get $u\in\C$.

It follows that $\C$ is the set of all regular elements of $A$, and each of 
these elements is invertible in $Q$. By Lemma 2.8 each right ideal in the set 
$\{\pi(I)\mid I\in\G_r\}$ contains a regular element of the ring $\pi(A)$. 
Therefore $I\cap\C\ne\varnothing$ for all $I\in\G_r$. Given any $q\in Q$, we 
have $qI\sbs A$ for some $I\in\G_r$, whence there exists $u\in\C$ such that 
$qu\in A$. This shows that $Q$ is a classical right quotient ring of $A$.

\medbreak
We have not verified yet that $Q$ is Artinian. For the proof of this fact we 
will need another consequence of assumption (A4):

\proclaim
Claim 4.
$\ \GKdim A/P=d(A)$ for each minimal prime ideal $P$ of $A$.
\endproclaim

Indeed, we have observed in the proof of Claim 2 that the prime radical $N$ of 
$A$ is the intersection of the prime ideals $P_i=M_i\cap A$. Therefore any 
minimal prime ideal $P$ of $A$ must coincide with one of these $P_i$. In other 
words, $P=M\cap A$ for some maximal ideal $M$ of $Q$. But then $A/P$ embeds in 
$Q/M\!$. Since $Q/M$ is $\G_r$-torsionfree by (A4), so too is $A/P$. This 
yields $\GKdim A/P\ge d(A)$, while the opposite inequality is clear from the 
definition of $d(A)$.

\medskip
We see that $\Om(A)$ is the set of all minimal prime ideals of $A$. 
Denoting by $Q(A/P)$ the classical quotient ring of the prime PI algebra 
$A/P$, we have
$$
Q/J\cong\prod_{P\in\Om(A)}Q(A/P)
$$
in view of Claim 3. Since $J$ is a two-sided ideal of $Q$ equal to $NQ$ by 
Claim 2, we have $J^i=N^iQ$, and since $Q$ is left flat over $A$ by a standard 
property of classical right quotient rings, there are isomorphisms 
$J^i\cong N^i\ot_AQ$ and
$$
J^i/J^{i+1}\cong N^i/N^{i+1}\ot_AQ/J
\cong\prod_{P\in\Om(A)}N^i/N^{i+1}\ot_AQ(A/P)
$$
for all $i\ge0$. The ideals $N^i$ form a Loewy series for the right $A$-module 
$A_A$. The finiteness of $\rho_P(A_A)$ implies that all right $Q/J$-modules 
$N^i/N^{i+1}\ot_AQ(A/P)$ have finite length. Hence so too do the $Q/J$-modules 
$J^i/J^{i+1}$, and we conclude that $Q$ has finite length as a right module 
over itself.

Thus $Q$ is right Artinian. By the left-right symmetry $Q$ is also left 
Artinian and is a classical left quotient ring of $A$.
\endproof

\section
3. Finitely generated PI module algebras

Let $H$ be a Hopf algebra over the base field $k$ and $A$ an \emph{$H$-module 
algebra}. This means that $A$ is a $k$-algebra equipped with a left $H$-module 
structure such that
$$
h(ab)=\sum\,(h\1a)\,(h\2b)\quad\hbox{for all $h\in H$ and $a,b\in A$},
$$
and $h1_A=\ep(h)1_A$ for all $h\in H$ where $\ep:H\to k$ is the counit.

The algebra $A$ is said to be \emph{$H$-semiprime} if $A$ has no nonzero 
nilpotent $H$-stable ideals, and $A$ is \emph{$H$-prime} if $A\ne0$ and 
$IJ\ne0$ for each pair of nonzero $H$-stable ideals $I$ and $J$ of $A$.

Our aim in this section is to prove that an $H$-prime $H$-module algebra $A$ 
has an Artinian classical quotient ring in the case when $A$ is a finitely 
generated PI algebra with the additional assumptions that $A$ satisfies the 
ACC on right and left annihilators and the action of $H$ on $A$ is 
\emph{locally finite} in the sense that $\,\dim Ha<\infty\,$ for all $a\in A$ 
(see Theorem 3.9). As an intermediate step we will investigate the symmetric 
quotient ring $Q=Q(A)$ introduced in section 2. We will verify that $A$ 
satisfies all conditions needed for an application of Proposition 2.9.

The antipode $S$ of $H$ is not assumed to be bijective. Because of this we 
need slightly more complicated constructions in which the Hopf subalgebra 
$S(H)$ of $H$ is taken into account. In some places $A$ is assumed to be 
$S(H)$-prime, which is the stronger requirement that $IJ\ne0$ for each pair of 
nonzero $S(H)$-stable ideals $I$ and $J$ of $A$. However, the difference with 
the $H$-primeness eventually disappears.

It has been observed in \cite{Sk99, Cor. 1.9} that each $S(H)$-submodule 
of a locally finite-dimensional $H$-module is automatically an $H$-submodule.  
This is an easy consequence of the fact that the antipode of the dual Hopf 
algebra $\Hd$ is always injective (see \cite{Sk06}). In particular, we have

\proclaim
Lemma 3.1.
If the action of $H$ on $A$ is locally finite{\rm,} then each $S(H)$-stable 
ideal of $A$ is $H$-stable{\rm,} and so the $H$-primeness of $A$ is equivalent 
to the $S(H)$-primeness.
\endproclaim

For an arbitrary algebra $A'$ and a coalgebra $C$ denote by $[C,A']$ the vector 
space $\Hom_k(C,A')$ equipped with the convolution multiplication. If 
$\dim C<\infty$, then $[C,A']\cong A'\ot C^*$ as algebras. Let $C\cop$ be $C$ 
with the opposite comultiplication. Then $[C\cop,A']$ is $\Hom_k(C,A')$ 
equipped with the multiplication
$$ 
(\xi\times\eta)(c)=\sum\xi(c\2)\eta(c\1),\qquad  
\xi,\eta\in\Hom_k(C,A'),\ \,c\in C.
$$
For an ideal $I$ of $A$ and a subcoalgebra $C$ of $H$ put
$$
I_C=\{a\in A\mid Ca\sbs I\},
$$
which is another ideal of $A$. Define a map $\th:A\to[C\cop,A/I]$, 
$\,a\mapsto\th_a$, setting
$$
\th_a(c)=S(c)a+I,\qquad a\in A,\ \,c\in C.
$$
We have $\th_{ab}=\th_a\times\th_b$ for $a,b\in A$ since 
$S(c)(ab)=\sum\,\bigl(S(c\2)a\bigr)\bigl(S(c\1)b\bigr)$ for all $c\in C$. Note 
that $\th_1$ is the identity element $c\mapsto\ep(c)+I$ of the algebra 
$[C\cop,A/I]$. Thus $\th$ is a homomorphism of unital algebras and
$$
\Ker\th=I_{S(C)}=\{a\in A\mid S(C)a\sbs I\}.
$$

Denote by $\F$ the set of all finite-dimensional subcoalgebras of $H$.

\proclaim
Lemma 3.2.
Let $C\in\F$. Then the algebra $B=[C\cop,A/I]$ is a finitely generated module 
over its subalgebra $\th(A)$ both on the left and on the right.
\endproclaim

\Proof.
We will prove that $B$ is generated as a $\th(A)$-module by the 
finite-dimensional subspace $C^*\sbs\Hom_k(C,A/I)$. Note that $\th$ factors as
$$
A\mapr{}[C\cop,A]\mapr{}[C\cop,A/I]
$$
where the first map is the special case of $\th$ obtained for the zero ideal 
of the algebra $A$, while the second map is induced by the canonical 
surjection $A\to A/I$. Thus it suffices to prove the lemma assuming that 
$I=0$. In this case we define a linear transformation $\Psi$ of the space 
$\Hom_k(C,A)$ by the rule
$$
(\Psi\xi)(c)=\sum S(c\2)\xi(c\1),\qquad\xi\in\Hom_k(C,A),\ \,c\in C.
$$
It is bijective with the inverse transformation
$(\Psi^{-1}\xi)(c)=\sum S^2(c\2)\xi(c\1)$.

The space $\Hom_k(C,A)$ is generated by $C^*$ as a left $A$-module with respect 
to the plain action of $A$ such that
$$
(a\xi)(c)=a\xi(c)\quad{\rm for}\ \,a\in A,\ \,\xi\in\Hom_k(C,A),\ \,c\in C.
$$
Since\quad$\displaystyle\Psi(a\xi)(c)=S(c\2)\bigl(a\xi(c\1)\bigr)
=\smash{\sum}\,\bigl(S(c\3)a\bigr)\bigl(S(c\2)\xi(c\1)\bigr)$
$$
\hphantom{Since\quad\Psi(a\xi)(c)=S(c\2)\bigl(a\xi(c\1)\bigr)\qquad}
{}=\smash{\sum}\,\th_a(c\2)(\Psi\xi)(c\1)=(\th_a\times\Psi\xi)(c)
$$
for all $c$, we have $\Psi(a\xi)=\th_a\times\Psi\xi$ for all $a$ and $\xi$. 
Thus $\Psi$ transforms the left module structure on $\Hom_k(C,A)$ given by the 
plain action of $A$ to the one obtained via the algebra homomorphism $\th$. 
Since $\Psi\xi=\xi$ for all $\xi\in C^*$, we get $B=\th(A)C^*$ in the algebra 
$B$. The other equality $B=C^*\mskip1mu\th(A)$ is obtained by means of similar 
arguments which use the bijective transformation $\Phi$ of the space 
$\Hom_k(C,A)$ defined by the rule $\,(\Phi\xi)(c)=\sum S(c\1)\xi(c\2)$.
\endproof

\proclaim
Lemma 3.3.
Let $A$ be a finitely generated $S(H)$-prime $H$-module PI algebra. If $P$ is 
a prime ideal of $A$ such that $\,\lann_A\mskip-2mu P\ne0,$ then $P_{S(C)}=0$ 
for some $C\in\F$. In this case $\,\th:A\to[C\cop,A/P]\,$ is injective.
\endproclaim

\Proof.
For each $C\in\F$ the left annihilator $\lann_AP_{S(C)}$ is a two-sided ideal 
of $A$ since so is $P_{S(C)}$. Denote by $V(C)$ the intersection of all prime 
ideals of $A$ containing $\lann_AP_{S(C)}$.

If $C,C'\in\F$ and $C\sbs C'$, then $P_{S(C)}\sps P_{S(C')}$, whence 
$V(C)\sbs V(C')$. Since $C_1+C_2\in\F$ for any two subcoalgebras 
$C_1,C_2\in\F$, the set of semiprime ideals $\{V(C)\mid C\in\F\}$ is directed 
by inclusion, i.e., for any two ideals in this set there is a larger ideal in 
the same set. Since $A$ satisfies the ACC on semiprime ideals by \cite{Pro, 
Ch. V, Cor. 2.2} or \cite{Row, Cor. 6.3.36$'$}, the set $\{V(C)\mid C\in\F\}$ 
has a largest element $V_0$. Thus $V(C)\sbs V_0$ for all $C\in\F$.

Put $K=\lann_A\mskip-2mu P$. If $a\in K$ and $b\in P_{S(C)}$, then
$$
(ca)b=\sum\,(c\1a)\bigl(c\2S(c\3)b\bigr)=\sum\,c\1\bigl(a\,S(c\2)b\bigr)=0
$$
for all $c\in C$ since $S(C)b\sbs P$. This shows that 
$CK\sbs\lann_AP_{S(C)}\sbs V(C)$. Hence $CK\sbs V_0$ for each $C\in\F$, and 
therefore $HK\sbs V_0$.

Denote by $I$ the ideal of $A$ generated by the $H$-submodule $HK$. Then $I$ 
is stable under the action of $H$ and $I\sbs V_0$. Since $K\ne0$ by the 
hypothesis, we also have $I\ne0$.

Now pick $C\in\F$ such that $V_0=V(C)$. Note that $V_0/\lann_AP_{S(C)}$ is the 
prime radical of the factor algebra $A/\lann_AP_{S(C)}$. Since the prime 
radicals of finitely generated PI algebras are nilpotent, we have 
$V_0^n\sbs\lann_AP_{S(C)}$ for some $n>0$. Then $I^n\sbs\lann_AP_{S(C)}$, 
and therefore $P_{S(C)}\sbs\rann_AI^n$.

Note that $I^n$ is an $H$-stable ideal of $A$. If $a\in I^n$ and 
$b\in\rann_AI^n$, then
$$
a\bigl(S(h)b\bigr)=\sum\,\bigl(S(h\2)h\3a\bigr)\,\bigl(S(h\1)b\bigr)
=\sum\,S(h\1)\bigl((h\2a)\,b\bigr)=0
$$
for all $h\in H$ since $Ha\sbs I^n$. This shows that the ideal $\rann_AI^n$ is 
stable under the action of $S(H)$. Since $A$ is $S(H)$-prime and $I\ne0$, we 
must have $\,\rann_AI^n=0$. Hence $P_{S(C)}=0$ too, and therefore $\Ker\th=0$.
\endproof

\Remark 2.
The conclusion of Lemma 3.3 remains true when $A$ satisfies the ACC on left 
annihilators, but is not necessarily finitely generated PI. This follows from 
the fact that the set of annihilators $\{\,\lann_AP_{S(C)}\mid C\in\F\}$ has a 
largest element, and this largest ideal in this set contains the $H$-stable 
ideal $I$ introduced in the proof of the lemma. By similar arguments it can be 
shown that Lemma 3.3 holds also in the case when $A$ satisfies the ACC on right 
annihilators and is $S^2(H)$-prime.
\endremark

\setitemsize(ii)
\proclaim
Lemma 3.4.
Let $A$ be a finitely generated $S(H)$-prime $H$-module PI algebra. Then\/{\rm:}

\item(i)
\ $\GKdim I=\GKdim A=d(A)$ for each nonzero one-sided ideal $I$ of $A,$

\item(ii)
\ $\rho_P(A_A)<\infty$ and $\rho_P({}_AA)<\infty$ for each $P\in\Om(A)$.

\endproclaim

\Proof.
Take any $P'\in\Spec A$ such that $\,\lann_A\mskip-2mu P'\ne0$. Such a prime 
ideal exists, as is seen from Lemma 1.3 applied to $M=A$. By Lemma 3.3 there 
is a coalgebra $C\in\F$ such that $A$ embeds in the algebra $B=[C\cop,A/P']$ 
by means of $\,\th:A\to B$.

One immediate consequence of this embedding is that $\GKdim A\le\GKdim B$. Now 
$B\cong A'\ot D$ where $A'=A/P'$ and $D$ is the algebra dual to the coalgebra 
$C\cop\!$. Since $\,\dim D<\infty$, we have
$$
\GKdim B=\GKdim A'\le\GKdim A.
$$
Thus an equality is attained here, which shows also that $\,\GKdim A=d(A)$.

Let $I$ be a right ideal of $A$. Since $B$ is a finitely generated left 
$\th(A)$-module according to Lemma 3.2, there is an inequality
$$
\GKdim_B(I\ot_AB)\le\GKdim_AI
$$
obtained by viewing $B$ as an $(A,B)$-bimodule and applying \cite{Mc-R, 8.3.14}. 
Since $\th(I)B$ is an epimorphic image of the right $B$-module $I\ot_AB$, we get
$$
\GKdim_B\bigl(\th(I)B\bigr)\le\GKdim_AI.
$$
We may identify $A'$ with the subalgebra $A'\ot1$ of $B\cong A'\ot D$. It is 
clear from the definition of Gelfand-Kirillov dimension that 
$\,\GKdim_BM\ge\GKdim_{A'}M\,$ for each right $B$-module $M$. The prime PI 
algebra $A'$ is GK-pure by \cite{Kr-L, Lemma 5.12}. It follows that 
$\,\GKdim_{A'}M=\GKdim A'\,$ whenever $M$ is a nonzero submodule of a free 
$A'$-module since such a module $M$ contains a submodule isomorphic to a nonzero 
right ideal of $A'$. In particular, this can be applied to right ideals of $B$ 
since $B$ is free as a right $A'$-module. If $I\ne0$, then $\th(I)B$ is a 
nonzero right ideal of $B$, and so
$$
\GKdim_B\bigl(\th(I)B\bigr)\ge\GKdim_{A'}\bigl(\th(I)B\bigr)=\GKdim A'.
$$
Hence $\,\GKdim_AI\ge\GKdim A'$. On the other hand, 
$$
\GKdim_AI\le\GKdim A=\GKdim A'.
$$
This shows that $\,\GKdim_AI=\GKdim A\,$. By the left-right symmetry such an 
equality holds also for each nonzero left ideal of $A$. This completes the 
proof of (i).

The prime PI algebra $A'$ has a simple Artinian classical quotient ring $Q'$ 
which is finite dimensional over its center $F$. By \cite{Kr-L, Th. 10.5} or 
\cite{Row, Prop. 6.3.40}
$$
\GKdim A'=\trdeg F/k.
$$
Hence $\trdeg F/k=d(A)$. Note that $F$ is a finitely generated extension field 
of $k$. This follows from the fact that there is an element $s\ne0$ in the 
center $Z$ of $A'$ such that the localization $A'[s^{-1}]$ is an Azumaya 
algebra over its center $Z[s^{-1}]$ (see \cite{Row, Cor. 6.1.36}. By the 
Artin-Tate Lemma $Z[s^{-1}]$ is a finitely generated algebra over $k$ since so 
is $A'[s^{-1}]$. Therefore the quotient field of $Z[s^{-1}]$, equal to $F$, 
is a finitely generated extension of $k$.

The algebra $B$, and therefore also the algebra $A$, embed in $Q'\ot D$. We 
may view $M=Q'\ot D$ as an $(F,A)$-bimodule. Since $\dim_FM<\infty$, Lemma 1.9 
can be applied. It yields $\rho_P(M)<\infty$ for each $P\in\Om(A)$. Since $A$ 
is an $A$-submodule of $M$, we get $\rho_P(A_A)<\infty$ for each $P\in\Om(A)$.  
Viewing $Q'\ot D$ as an $(A,F)$-bimodule, we deduce similarly the other part 
of (ii).
\endproof

The action of $H$ on $A$ gives rise to certain tensoring operations on the 
categories of right and left $A$-modules. Our subsequent arguments will use 
the fact that the classes of $\G_r$-torsion and $\G_l$-torsion modules are 
stable under those operations.

If $M$ is a right $A$-module and $U$ is a left $H$-comodule, then the vector 
space $M\mskip-1mu\ot U$ is a right $A$-module with respect to the 
\emph{twisted action} of $A$ defined by the rule
$$
(m\ot u)a=\sum m(u_{(-1)}a)\ot u\0\mskip 1mu,\qquad m\in M,\ u\in U,\ a\in A.
$$
For a left $A$-module $M$ and a right $H$-comodule $U$ there is a twisted left 
$A$-module structure on the vector space $M\ot U$ defined by the rule
$$
a(m\ot u)=\sum\,(u\1a)m\ot u\0\mskip 1mu.
$$
Here the symbolic notation $\sum u_{(-1)}\ot u\0\in H\ot U$ and 
$\sum u\0\ot u\1\in U\ot H$ is used to represent the values of the comodule 
structure maps $U\to H\ot U$ and $U\to U\ot H$, respectively. Note that the 
above tensoring operations make sense for an arbitrary $H$-module algebra $A$.

\proclaim
Lemma 3.5.
Let $M$ be a right $A$-module and $U$ a left $H$-comodule. Then
$$
\GKdim M\ot U\le\GKdim M.
$$
Suppose that $A$ is a finitely generated $H$-module PI algebra. Then 
$M\ot U\in\T_r$ whenever $M\in\T_r$.

Similar conclusions hold for a left $A$-module and a right $H$-comodule.
\endproclaim

\Proof.
Let $V\sbs A$ and $X\sbs M\ot U$ be finite-dimensional subspaces with $1\in V$. 
We can find a subspace $Y\sbs M$ and a subcomodule $U'\sbs U$, both of finite 
dimension, such that $X\sbs Y\ot U'$. There is a coalgebra $C\in\F$ such that
$$
\sum u_{(-1)}\ot u\0\in C\ot U'
$$
for all $u\in U'$. Put $W=CV$, which is a finite-dimensional subspace of $A$. 
Note that $CV^n\sbs W^n$ for each integer $n>0$. Hence
$$
XV^n\sbs(Y\ot U')V^n\sbs Y(CV^n)\ot U'\sbs YW^n\ot U',
$$
which yields $\,\dim XV^n\le(\dim U')(\dim YW^n)\,$ and
$$
\limsup_{n\to\infty}\,\log_n(\dim XV^n)\le
\limsup_{n\to\infty}\,\log_n(\dim YW^n)\le\GKdim M.
$$
Taking the supremum over all $V$ and $X$, we get the required inequality.

If $M\in\T_r,$ then $M$ has a chain of submodules 
$M=M_0\sps M_1\sps\ldots\sps M_n=0$ such that $\GKdim M_{i-1}/M_i<d(A)$ for 
all $i$. The spaces $M_i\ot U$ form a chain of submodules in $M\ot U$. 
For each $i$ we have
$$
(M_{i-1}\ot U)/(M_i\ot U)\cong(M_{i-1}/M_i)\ot U,
$$
which is an $A$-module of Gelfand-Kirillov dimension less than $d(A)$ by 
the already established inequality. It follows that $M\ot U\in\T_r$.
\endproof

\proclaim
Corollary 3.6.
The action of $H$ on $A$ is $\G_r$-continuous and $\G_l$-continuous in the 
sense that for every $h\in H,$ $I\in\G_r,$ $J\in\G_l$ there exist $I'\in\G_r$ 
and $J'\in\G_l$ such that $hI'\sbs I$ and $hJ'\sbs J$.
\endproclaim

\Proof.
The right $A$-module $M=A/I$ is $\G_r$-torsion since $I\in\G_r$. Denote by 
$v$ the coset $1+I\in M$. Consider $H$ as a left $H$-comodule with 
respect to the comultiplication in $H$. Then $M\ot H$ is a right $A$-module 
with respect to the twisted action of $A$. It is $\G_r$-torsion by Lemma 3.5. 
Hence the element $v\ot h\in M\ot H$ is annihilated by some right 
ideal $I'\in\G_r$. If $a\in I'$, then
$$
\sum\,v(h\1a)\ot h\2=(v\ot h)a=0\quad\hbox{in $M\ot H$},
$$
and it follows that $v(ha)=\sum\ep(h\2)\,v(h\1a)=0$. Since 
$I$ is the annihilator of $v$ in $A$, we see that $ha\in I$ for all 
$a\in I'$, as required.

The other conclusion is obtained in a similar way by considering the left 
$A$-module $M\ot H$ where $M=A/J$ and $H$ is regarded as a right $H$-comodule.  
\endproof

\proclaim
Proposition 3.7.
Let $A$ be a finitely generated $S(H)$-prime $H$-module PI algebra satisfying 
the ACC on right and left annihilators.

Then its symmetric quotient ring $Q$ with respect to the Gabriel topologies 
$\G_l,\mskip1mu\G_r$ is a semiprimary $S(H)$-prime $S(H)$-module algebra.
\endproclaim

\Proof.
Assumptions (A2), (A3) of section 2 are satisfied by Lemma 3.4, while (A1) is 
included in the hypothesis. By Lemma 2.6 $Q$ is semiprimary.

Since the action of $H$ on $A$ is $\G_r$-continuous, it extends to the right 
quotient ring $Q_r$ of $A$ with respect to $\G_r$, as explained in 
\cite{Mo-Sch95, Th. 3.13}. Explicitly, let $h\in H$, and let $q\in Q_r$ be 
represented by a right $A$-linear map $f:I\to A$ where $I\in\G_r$. There is a 
coalgebra $C\in\F$ containing $h$. Since $\dim C<\infty$, we can use Corollary 
3.6 to find a right ideal $I'\in\G_r$ such that $S(c)I'\sbs I$ for all 
$c\in C$. Then $hq\in Q_r$ is represented by the right $A$-linear map $I'\to A$ 
given by the assignment
$$
a\mapsto\sum h\1f\bigl(S(h\2)a\bigr),\qquad a\in I'.
$$
This action makes $Q_r$ into an $H$-module algebra containing $A$ as an 
$H$-stable subalgebra. (Such an action was first introduced by Cohen 
\cite{Coh85, Th. 18} in the case of the quotient ring constructed with respect 
to the filter of $H$-stable two-sided ideals with zero left and right 
annihilators.)

If $q\in Q$, then $Jq\sbs A$ for some $J\in\G_l$. Given $h$ and $C$ as above, 
there exists $J'\in\G_l$ such that $cJ'\sbs J$ for all $c\in C$, again by 
Corollary 3.6. Then
$$
a\bigl(S(h)q\bigr)
=\sum\,\bigl(S(h\2)h\3a\bigr)\,\bigl(S(h\1)q\bigr)
=\sum\,S(h\1)\bigl((h\2a)\,q\bigr)\in A
$$
for all $a\in J'$. Hence $J'\,S(h)q\sbs A$, and therefore $S(h)q\in Q$. This 
shows that $Q$ is an $S(H)$-stable subalgebra of $Q_r$.

If $I$ is any nonzero $S(H)$-stable ideal of $Q$, then $I\cap A$ is a nonzero 
$S(H)$-stable ideal of $A$. The condition that $A$ is $S(H)$-prime implies that 
$(I\cap A)(J\cap A)\ne0$, and therefore $IJ\ne0$, for each pair of nonzero 
$S(H)$-stable ideals of $Q$. Hence $Q$ is $S(H)$-prime.
\endproof

Proposition 3.7 does not give us all we want. It is not clear whether one can 
achieve more under the same assumptions about $A$. At this point we are going 
to employ the local finiteness of the action.

Recall that an $H$-module algebra is said to be \emph{$H$-simple} if it has 
no $H$-stable ideals other than the zero ideal and the whole algebra. Our 
further arguments are based on the already known result for which we cite 
\cite{Sk99, Prop. 3.2}:

\proclaim
Proposition 3.8.
Suppose that $Q$ is a semiprimary $H$-semiprime $H$-module algebra containing 
an $H$-stable subalgebra $A$ such that the action of $H$ on $A$ is locally 
finite and $I\cap A\ne0$ for each nonzero right ideal $I$ of $Q$. Then
there is an isomorphism of $H$-module algebras
$$
Q\cong Q_1\times\cdots\times Q_n
$$
where $Q_1,\ldots,Q_n$ are $H$-simple $H$-module algebras.
\endproclaim

\setitemsize(iii)
\proclaim
Theorem 3.9.
Let $A$ be a finitely generated $H$-prime $H$-module PI algebra satisfying 
the ACC on right and left annihilators. Suppose that the action of $H$ on $A$ 
is locally finite. Then\/{\rm:}

\item(i)
\ $\GKdim I=\GKdim A$ for each nonzero one-sided ideal $I$ of $A,$

\item(ii)
\ $\GKdim A/P=\GKdim A$ for each minimal prime ideal $P$ of $A,$

\item(iii)
\ $A$ has a right and left Artinian classical quotient ring $Q,$

\item(iv)
\ $Q$ is an $H$-simple $H$-module algebra.

\endproclaim

\Proof.
By Lemma 3.1 $A$ is $S(H)$-prime. Hence Lemma 3.4 and Proposition 3.7 apply to 
$A$. In particular, assumptions (A1)--(A3) of section 2 are satisfied, and (i) 
holds by Lemma 3.4. Let $Q$ be the symmetric quotient ring of $A$ constructed 
in section 2. We will verify condition (A4) in the statement of Proposition 
2.9. This will establish assertion (iii), and (ii) will follow at once from 
Lemma 3.4 and Claim 4 in the proof of Proposition 2.9.

By Proposition 3.7 $Q$ is a semiprimary $S(H)$-prime $S(H)$-module algebra. 
The algebra $A$ embeds in $Q$ as an $S(H)$-stable subalgebra. By the 
construction of $Q$ each nonzero right $A$-submodule of $Q$ has nonzero 
intersection with $A$. Hence we can apply Proposition 3.8 with $H$ replaced by 
$S(H)$. It follows that $Q\cong Q_1\times\cdots\times Q_n$ where each $Q_i$ is 
an $S(H)$-simple $S(H)$-module algebra. Since $Q$ is $S(H)$-prime, we must 
have $n=1$, and so $Q$ is $S(H)$-simple.

For a right $A$-module $M$ and a left $H$-comodule $U$ we have defined the 
twisted action of $A$ on the vector space $M\ot U$. If $M$ is a right $Q$-module 
and $U$ is a left $S(H)$-comodule, then $M\ot U$ is a right $Q$-module with 
respect to the twisted action of $Q$. Take $U=S(H)$, which is a left 
$S(H)$-comodule with respect to the comultiplication. If $I$ is the 
annihilator of $M$ in $Q$, then the annihilator of $M\ot S(H)$ in $Q$ is the 
$S(H)$-stable ideal
$$
I_{S(H)}=\{q\in Q\mid S(H)\mskip1mu q\sbs I\}.
$$
Indeed, if $q\in Q$ annihilates $M\ot S(H)$, then for all $m\in M$ and 
$h\in S(H)$ we have
$$
0=(m\ot h)q=\sum m(h\1q)\ot h\2\quad\hbox{in $M\ot S(H)$},
$$
and applying the map $\id\ot\ep$, we get $m(hq)=0$, which shows that 
$q\in I_{S(H)}$. Conversely, it is seen from the formula for the twisted 
action that all elements of $I_{S(H)}$ annihilate $M\ot S(H)$.

If $M\ne0$, then $I_{S(H)}=0$ since $Q$ is $S(H)$-simple. In this case the 
$Q$-module $M\ot S(H)$ is faithful. Suppose that $M$ is $\G_r$-torsion. Then 
so is $M\ot S(H)$ by Lemma 3.5. Since $Q$ is semiprimary, it satisfies the DCC 
on finitely generated right ideals \cite{St, Ch. VIII, Prop. 5.5}. Hence $Q$ has 
a simple right ideal, say $V$. Since $V$ does not annihilate $M\ot S(H)$, it 
embeds in $M\ot S(H)$ as a $Q$-submodule. But then $V$ has to be $\G_r$-torsion, 
which contradicts the fact that $Q$ is $\G_r$-torsionfree (see Lemma 2.1).

The preceding argument shows that there exist no nonzero right $Q$-modules which 
are $\G_r$-torsion. Each right $A$-module $M$ has a largest $\G_r$-torsion 
submodule $T$, and the factor module $M/T$ is $\G_r$-torsionfree. Suppose that 
$M$ is a right $Q$-module. Given $q\in Q$, there is $I\in\G_r$ such that 
$qI\sbs A$. Since then $Tq\mskip1mu I\sbs TA\sbs T$, we must have $Tq\sbs T$ 
by $\G_r$-torsionfreeness of $M/T$. Therefore $T$ is a $Q$-submodule of $M$, 
which is possible only when $T=0$.

We have thus proved that each right $Q$-module is $\G_r$-torsionfree. Using 
tensoring operations on left $Q$-modules, it is proved in a similar way that 
each left $Q$-module is $\G_l$-torsionfree.

Now we can use Proposition 2.9 and all parts in its proof. There we have seen 
that $Q$ is a classical quotient ring of $A$ and each right ideal $I\in\G_r$ 
contains a regular element of $A$. This implies that $Q=Q_r$. Indeed, for each 
$q\in Q_r$ there exists a regular element $s\in A$ such that $qs\in A$, but 
then $q\in Q$ since $s^{-1}\in Q$.

In the proof of Proposition 3.7 we have seen that the action of $H$ on $A$ 
extends to $Q_r$. Therefore $Q$ is not only an $S(H)$-module algebra, but an 
$H$-module algebra. It is $H$-simple by Proposition 3.8. This proves the 
remaining assertion (iv).
\endproof

The existence of an Artinian classical quotient ring proved in Theorem 3.9 
for an $H$-prime algebra can be extended to the $H$-semiprime case. For the 
sake of completeness we present this result below, although it will not be 
needed later in this paper. Such an extension makes use of two standard 
ring-theoretic facts.

\proclaim
Lemma 3.10.
Let $I_1,\ldots,I_n$ be ideals of a ring $R$ such that $I_1\cap\cdots\cap I_n=0$ 
and for each $i$ the factor ring $R/I_i$ has a classical {\rm(}left{\rm,} 
right{\rm)} quotient ring $Q_i$. Suppose that $I_jQ_i=Q_i$ for each pair of 
indices $i\ne j$. Then $Q=Q_1\times\cdots\times Q_n$ is a classical 
{\rm(}left{\rm,} right{\rm)} quotient ring of $R$.
\endproclaim

This conclusion is better known in the case when all factor rings $R/I_i$ are 
prime Goldie. However, the proof of \cite{Row, Th. 3.2.27} can be easily 
adapted for the more general version stated in Lemma 3.10 above. The equality 
$I_jQ_i=Q_i$ means that $I_j$ contains an element whose image in $R/I_i$ is a 
regular element of that ring. Therefore, for each $i$ it is possible to find 
$a_i\in\bigcap_{j\ne i}I_j$ such that $a_i+I_i$ is regular in $R/I_i$. Then 
the arguments given in \cite{Row} show that each regular element of $R$ is 
invertible in $Q$, and each element $q\in Q$ can be written as $b^{-1}r$ 
(respectively, $rb^{-1}$) with $b,r\in R$ in the case of left (respectively, 
right) quotient rings.

\proclaim
Lemma 3.11.
Let $L$ be the left annihilator of a two-sided ideal of a ring $R$. If $R$ 
satisfies the ACC on right and left annihilators{\rm}, then so does the factor 
ring $R/L$.
\endproclaim

For this conclusion we cite \cite{Her-S64, Lemma 3}.

\proclaim
Theorem 3.12.
Let $A$ be a finitely generated $H$-semiprime $H$-module PI algebra satisfying 
the ACC on right and left annihilators. Suppose that the action of $H$ on $A$ 
is locally finite. Then $A$ has a right and left Artinian classical quotient 
ring.
\endproclaim

\Proof.
Let $P_1,\ldots,P_m$ be all the minimal prime ideals of $A$. Their intersection 
is the prime radical $N$ of $A$. Recall that $N$ is nilpotent. Each $P_i$ 
contains an $H$-prime $H$-stable ideal
$$
(P_i)_H=\{a\in A\mid Ha\sbs P_i\}.
$$
By $H$-semiprimeness of $A$ we have $(P_1)_H\cap\cdots\cap(P_m)_H=0$ since 
this intersection is an $H$-stable ideal contained in $N$.

Let $I_1,\ldots,I_n$ be all the minimal ideals among $(P_1)_H,\ldots,(P_m)_H$. 
Then each $I_i$ is $H$-prime, and $I_1\cap\cdots\cap I_n=0$. Put 
$K_i=\bigcap_{j\ne i}I_j$. Since $K_i$ is an $H$-stable ideal of $A$, so is 
its left annihilator $L_i$. Clearly $I_i\sbs L_i$. Since $I_i$ is $H$-prime 
and $I_j\not\sbs I_i$ whenever $j\ne i$, we have $K_i\not\sbs I_i$. Then the 
inclusion $L_iK_i=0\sbs I_i$ entails $L_i\sbs I_i$. It follows that 
$I_i=\lann_AK_i$.

For each $i$ the finitely generated $H$-prime $H$-module PI algebra $A/I_i$ 
satisfies the ACC on right and left annihilators in view of Lemma 3.11. Hence 
$A/I_i$ has a right and left Artinian classical quotient ring $Q_i$ by Theorem 
3.9. Also, $Q_i$ is an $H$-simple $H$-module algebra. If $j\ne i$, then 
$I_jQ_i\ne0$ since $I_j\not\sbs I_i$. On the other hand, $I_jQ_i$ is a 
two-sided ideal of $Q_i$ (see \cite{Mc-R, 2.1.16}). Since this ideal is stable 
under the action of $H$, we must have $I_jQ_i=Q_i$. Now Lemma 3.10 shows that
$Q_1\times\cdots\times Q_n$ is a classical quotient ring of $A$.
\endproof

\section
4. Flatness over coideal subalgebras and Hopf subalgebras

Let $H$ be a Hopf algebra over the field $k$ with the comultiplication $\De$, 
the counit $\ep$ and the antipode $S$. Any subalgebra $A\sbs H$ such that 
$\De(A)\sbs A\ot H$ is called a \emph{right coideal subalgebra}.

If $H$ has a right Artinian classical right quotient ring $Q(H)$, then $S$ is 
bijective by \cite{Sk06, Th. A}. In this case $S$ is an antiautomorphism of 
$H$, and $S$ extends to an antiautomorphism of $Q(H)$, which implies that 
$Q(H)$ is also a left quotient ring and is left Artinian. Therefore it will 
not lessen generality of results if we use two-sided conditions on $Q(H)$ 
instead of one-sided. By an \emph{Artinian ring} we mean a ring which is left 
and right Artinian.

\setitemsize(iii)
\proclaim
Theorem 4.1.
Let $A$ be a PI right coideal subalgebra of a residually finite-dimen\-sional 
Hopf algebra $H$ which has an Artinian classical quotient ring $Q(H)$. Then\/{\rm:}

\item(i)
\ $A$ has an Artinian classical quotient ring $Q(A),$

\item(ii)
\ $H$ is left and right $A$-flat{\rm,}

\item(iii)
\ if $A$ is a Hopf subalgebra{\rm,} then $H$ is left and right faithfully $A$-flat.

\endproclaim

\Proof.
Let $A\op$ and $H\op$ be $A$ and $H$ with the opposite multiplications. Note 
that $H\op$ is a residually finite-dimensional Hopf algebra, and $A\op$ a PI 
right coideal subalgebra of $H\op$. Passage from $A$ to $A\op$ shows that the 
left versions of the properties in (i)---(iii) are equivalent to their right 
versions. It is now clear that (i)$\Rar$(ii) by \cite{Sk10, Th. 1.8} or 
\cite{Sk99, Th. 4.4}. Also, (ii)$\Rar$(iii) by \cite{Ma-W94, Th. 2.1}. 
(Moreover, (iii) implies that $H$ is a projective generator as either right or 
left $A$-module \cite{Schn92, Cor. 1.8}.)

Thus only (i) has to be proved. Suppose first that $A$ is a finitely generated 
algebra. Consider the dual Hopf algebra $\Hd$ which consists of all linear 
functions $H\to k$ vanishing on an ideal of finite codimension in $H$. There 
is a locally finite action of $\Hd$ on $A$ defined by the rule
$$
f\rhu a=\sum\,f(a\2)\,a\1\,,\qquad f\in\Hd,\ \,a\in A.
$$
With respect to this action $A$ is an $\Hd$-prime $\Hd$-module algebra (see 
\cite{Sk99, Lemma 4.3} for an easy verification). Since $A$ embeds in the 
Artinian ring $Q(H)$, it satisfies the ACC on right and left annihilators. 
Therefore (i) follows from Theorem 3.9.

Consider now the case when $A$ is not finitely generated. Denote by $\I$ the 
set of all finitely generated right coideal subalgebras of $H$ contained in 
$A$. The subalgebras generated by finite-dimensional right $H$-subcomodules of 
$A$ belong to $\I$. Since each element of $A$ is contained in a 
finite-dimensional right subcomodule, the set $\I$ is directed by inclusion 
and the union of all subalgebras $A'\in\I$ give the whole $A$. It follows that
$$
V\ot_AH\cong\limdir_{A'\in\I}V\ot_{A'}H
$$
for each right $A$-module $V$. We know already that $H$ is left $A'$-flat for 
each $A'\in\I$. Hence the functors $?\ot_{A'}H$ are exact, and so is $?\ot_AH$ 
by exactness of inductive direct limits. In other words, $H$ is left $A$-flat. 
By symmetry $H$ is right $A$-flat. This proves (ii), and (iii) also follows.

To complete the proof of (i) we need some facts which will be stated separately. 
Each algebra $A'\in\I$ has an Artinian quotient ring $Q(A')$, as we have 
established already. By Lemma 4.3 $Q(A')$ embeds in $Q(H)$ as a subring. Put 
$$
Q(A)=\bigcup_{A'\in\I}Q(A')\sbs Q(H).
$$
If $A'\in\I$, then each regular element of $A'$ is invertible in $Q(H)$, and 
therefore regular in $A$. Conversely, each regular element $s$ of $A$ is 
regular in any subalgebra $A'$ containing $s$. Hence all regular elements of 
$A$ are invertible in $Q(H)$. If $q\in Q(A)$, then $q\in Q(A')$ for some 
$A'\in\I$, whence $qs\in A'\sbs A$ for some regular element $s\in A'$. This 
shows that $Q(A)$ is a classical right quotient ring of $A$. By symmetry 
$Q(A)$ is also a classical left quotient ring of $A$.

By Lemma 4.4 $Q(H)$ is faithfully flat as a left $Q(A')$-module for any 
$A'\in\I$. Therefore the map $I\mapsto I\mskip1mu Q(H)$ embeds the lattice 
$\L(A')$ of right ideals of $Q(A')$ into the lattice $\L(H)$ of right ideals 
of $Q(H)$. Denote by $\L_f(A)$ the lattice of finitely generated right ideals 
of $Q(A)$. Given $I\in\L_f(A)$, there exists a subalgebra $A'\in\I$ 
containing a finite set of generators of $I$. Then $I=I'\mskip1mu Q(A)$ and 
$I\mskip1mu Q(H)=I'\mskip1mu Q(H)$ for some $I'\in\L(A')$. For two right ideals 
$I,J\in\L_f(A)$ we can find $A'\in\I$ such that both $I$ and $J$ are extensions 
of right ideals of $Q(A')$. It follows that $\L_f(A)$ embeds in $\L(H)$. Since 
$Q(H)$ is right Artinian, we deduce that $\L_f(A)$ satisfies the ACC and the DCC.

By the ACC on finitely generated right ideals an arbitrary right ideal $I$ of 
$Q(A)$ must coincide with the right ideal generated by a finite subset of $I$. 
Thus the lattice $\L_f(A)$ contains all right ideals of $Q(A)$, and the DCC 
for this lattice implies that $Q(A)$ is right Artinian. By symmetry $Q(A)$ is 
left Artinian.
\endproof

The next three lemmas provide additions needed to complete the proof of part 
(i) in Theorem 4.1.

\proclaim
Lemma 4.2.
Let $H$ be an arbitrary Hopf algebra{\rm,} $D$ an $H$-module algebra{\rm,} and
$B$ an $H$-stable subalgebra of $D$. Suppose that $B$ is semiprimary and 
$S(H)$-simple. If $D$ is left $B$-flat{\rm,} then $D$ is left faithfully 
$B$-flat.
\endproclaim

\Proof.
Denote by $\T_{D/B}$ the class of all right $B$-modules $V$ such that 
$V\ot_BD=0$. We have to show that $\T_{D/B}$ contains no nonzero modules. If 
$V\in\T_{D/B}$, then $\T_{D/B}$ contains all submodules of $V$ by the flatness 
assumption. (In fact $\T_{D/B}$ is the torsion class of a hereditary torsion 
theory.) The arguments that follow are similar to the proof that there exist 
no nonzero $\G_r$-torsion $Q$-modules in the setup of Theorem 3.9.

First we show that $\T_{D/B}$ is closed under tensoring with right comodules 
(note a difference with section 3 where right modules were tensored 
with left comodules). For a right $H$-comodule $U$ and a right $B$-module $V$ 
we will consider the vector space $U\ot V$ as a right $B$-module with respect 
to the twisted action of $B$ defined as
$$
(u\ot v)b=\sum u\0\ot v\bigl(S(u\1)b\bigr),\qquad u\in U,\ v\in V,\ b\in B.
$$
Similar tensoring operations are defined for right $D$-modules. With these 
conventions there is an isomorphism of right $D$-modules
$$
U\ot(V\ot_BD)\cong(U\ot V)\ot_BD
$$
constructed as follows. We may identify $V\ot_BD$ with the factor space of 
$V\ot D$ by its subspace spanned by $\{v\ot bd-vb\ot d\mid 
\hbox{$v\in V$, $b\in B$, $d\in D$}\}$. Then
$$
U\ot(V\ot_BD)\cong(U\ot V\ot D)/R_1\,,\qquad
(U\ot V)\ot_BD\cong(U\ot V\ot D)/R_2
$$
where $R_1$ and $R_2$ are the subspaces of $U\ot V\ot D$ spanned, 
respectively, by
$$
\openup1\jot
\eqalign{
&\{u\ot v\ot bd-u\ot vb\ot d\mid
\hbox{$u\in U$, $v\in V$, $b\in B$, $d\in D$}\}\qquad{\rm and}\cr
&\{u\ot v\ot bd-{\textstyle\sum}\,u\0\ot v\bigl(S(u\1)b\bigr)\ot d\mid
\hbox{$u\in U$, $v\in V$, $b\in B$, $d\in D$}\}.\cr
}
$$
Define linear transformations $\Phi$ and $\Psi$ of $U\ot V\ot D$ by the rules
$$
\Phi(u\ot v\ot d)=\sum u\0\ot v\ot u\1d\,,\qquad
\Psi(u\ot v\ot d)=\sum u\0\ot v\ot S(u\1)d.
$$
Then $\Psi=\Phi^{-1}$. For $u\in U$, $v\in V$, $b\in B$, $d\in D$ we have
$$
\openup1\jot
\eqalign{
\Phi(u\ot v\ot bd)&{}=\sum u\0\ot v\ot(u\1b)(u\2d)\cr
&{}\equiv\sum u\0\ot v\bigl(S(u\1)u\2b\bigr)\ot u\3d\cr
&\qquad\qquad{}=\sum u\0\ot vb\ot u\1d=\Phi(u\ot vb\ot d)\quad
\hbox{modulo $R_2$}\,,\cr
}
$$
$$
\openup1\jot
\eqalign{
\Psi(u\ot v\ot bd)&{}=\sum u\0\ot v\ot\bigl(S(u\2)b\bigr)\bigl(S(u\1)d\bigr)\cr
&{}\equiv\sum u\0\ot v\bigl(S(u\2)b\bigr)\ot S(u\1)d
=\sum\Psi\bigl(u\0\ot v\bigl(S(u\1)b\bigr)\ot d\bigr)\cr
}
$$
modulo $R_1$. It follows that $\Phi(R_1)\sbs R_2$ and $\Psi(R_2)\sbs R_1$. 
Hence $\Phi(R_1)=R_2$. Since $\Phi$ commutes with the action of $D$ by right 
multiplications on the third tensorand, we see that $\Phi$ induces the 
required isomorphism of $D$-modules.

If $V\ot_BD=0$, then $(U\ot V)\ot_BD=0$. In other words, $U\ot V\in\T_{D/B}$ 
whenever $V\in\T_{D/B}$, as claimed.

Now take $U=H$ with the comodule structure given by the comultiplication. The 
annihilator of the twisted right $B$-module $H\ot V$ is the largest 
$S(H)$-stable ideal of $B$ contained in the annihilator of $V$. Suppose that 
$V\ne0$. Then this ideal is zero since $B$ is $S(H)$-simple. In other words, 
$H\ot V$ is a faithful $B$-module.

Since $B$ is semiprimary, it has a simple right ideal, say $K$. This right 
ideal embeds in $H\ot V$ as a submodule. If $V\in\T_{D/B}$, then $H\ot 
V\in\T_{D/B}$, whence $K\in\T_{D/B}$ as well. But this is impossible since 
$K\ot_BD\cong KD\sps K$, and so $K\ot_BD\ne0$.

It follows that $V\notin\T_{D/B}$, i.e., $V\ot_BD\ne0$, for each nonzero right 
$B$-module $V$. Thus the flat left $B$-module $D$ is faithfully flat.
\endproof

\proclaim
Lemma 4.3.
Let $A$ be a right coideal subalgebra of a residually finite-dimensional Hopf 
algebra $H$. Suppose that $A$ and $H$ have right Artinian classical right 
quotient rings $Q(A)$ and $Q(H)$. Then the inclusion $A\hrar H$ extends to an 
injective homomorphism of $\Hd$-simple $\Hd$-module algebras $Q(A)\to Q(H)$.
\endproclaim

\Proof.
First, the action of $\Hd$ on $A$ and $H$ extends to $Q(A)$ and $Q(H)$ by 
\cite{Sk-Oy06, Th. 2.2}. The $\Hd$-module algebras $Q(A)$ and $Q(H)$ are 
$\Hd$-prime since so are $A$ and $H$, and by Proposition 3.8 they are in fact 
$\Hd$-simple (see also Claim 1 in the proof of \cite{Sk99, Th. 4.4}). The 
conclusion of Lemma 4.3 can now be deduced from \cite{Sk10, Lemma 1.7}, but 
the proof in this special case is easier.

To obtain a ring homomorphism $Q(A)\to Q(H)$ we have to show that all regular 
elements of $A$ remain regular in $H$. Denote by $I$ the set of all elements 
of $H$ with the property that the right annihilator of the element contains a 
regular element of $A$. Then $I$ is the kernel of the linear map $f:H\to 
H\ot_AQ(A)$ such that $f(x)=x\ot1$ for $x\in H$. Consider $H\ot Q(A)$ as a 
left $\Hd$-module with respect to the tensor product of the $\Hd$-module 
structures on $H$ and $Q(A)$. Then $H\ot_AQ(A)$ is its factor module, and the 
map $f$ is $\Hd$-linear. Hence $I$ is an $\Hd$-submodule of $H$, which means 
precisely that $\De(I)\sbs I\ot H$, i.e., $I$ is a right coideal of $H$. 
Obviously, $I$ is also a left ideal of $H$. Since the antipode of $H$ is 
bijective, we conclude that $I=0$. Indeed,
$$
\ep(x)1=\sum S^{-1}(x\2)x\1\in HI\sbs I
$$
for each $x\in I$. However, $1\notin I$ since the identity element has zero 
right annihilator. Hence $\ep(I)=0$, but then $x=\sum\ep(x\1)x\2=0$ for each 
$x\in I$.

Let $s$ be a regular element of $A$. The equality $I=0$ forces $\,\lann_Hs=0$. 
It follows that $\,\lann_{Q(H)}s$ is a left ideal of $Q(H)$ which has zero 
intersection with $H$. Since $Q(H)$ is a left quotient ring of $H$, this 
entails $\,\lann_{Q(H)}s=0$, i.e., $s$ is left regular in $Q(H)$. Using the 
fact that all left regular elements of any Artinian ring are invertible, we 
deduce that each regular element of $A$ is invertible in $Q(H)$.

By the universality property of classical quotient rings there is a unique 
ring homomorphism $\psi:Q(A)\to Q(H)$ extending the inclusion $A\hrar H$. The 
fact that $\psi$ is $\Hd$-linear is proved as in \cite{Sk10, Lemma 1.7}. Since 
$\Ker\psi$ has zero intersection with $A$, we get $\Ker\psi=0$. Hence $\psi$ 
is injective.
\endproof

\proclaim
Lemma 4.4.
Under the assumptions of Lemma\/ {\rm4.3} $Q(H)$ is faithfully flat as a left 
$Q(A)$-module.
\endproclaim

\Proof.
By Lemma 4.3 $Q(A)$ is identified with an $\Hd$-stable subalgebra of $Q(H)$, 
and this subalgebra is $\Hd$-simple. By \cite{Sk10, Th. 1.8} $H$ is left 
$A$-flat. Since $Q(H)$ is left $H$-flat, $Q(H)$ is also left $A$-flat by 
transitivity of flatness. Noting that
$$
V\ot_{Q(A)}Q(H)\cong\bigl(V\ot_AQ(A)\bigr)\ot_{Q(A)}Q(H)\cong V\ot_AQ(H)
$$
for each right $Q(A)$-module $V$, we deduce that $Q(H)$ is left $Q(A)$-flat. 
As we have noted earlier, the existence of a right Artinian classical right 
quotient ring $Q(H)$ implies that the antipode of $H$ is bijective. Then so is 
the antipode of $\Hd$. Now Lemma 4.2 applies with $B=Q(A)$, $D=Q(H)$, and $H$ 
replaced by $\Hd$.
\endproof

If $H$ is a residually finite-dimensional Noetherian Hopf algebra, then the 
existence of an Artinian classical quotient ring $Q(H)$ has been established 
in \cite{Sk99, Th. 4.5}. Thus Theorem 0.1 stated in the introduction is a 
special case of Theorem 4.1. And here is another special case where $H$ is not 
necessarily Noetherian:

\proclaim
Theorem 4.5.
Let $H$ be a residually finite-dimensional PI Hopf algebra{\rm,} finitely 
generated as an ordinary algebra. Suppose that $H$ satisfies the ACC on right 
and left annihilators. Then\/{\rm:}

\item(i)
\ $H$ has an Artinian classical quotient ring{\rm,}

\item(ii)
\ the antipode $S:H\to H$ is bijective{\rm,}

\item(iii)
\ $H$ is flat as either right or left module over any right coideal subalgebra{\rm,}

\item(iv)
\ $H$ is faithfully flat as either right or left module over any Hopf subalgebra.

\endproclaim

\Proof.
Since $H$ is an $\Hd$-prime (actually $\Hd$-simple) $\Hd$-module algebra, we 
can apply Theorem 3.9 to get (i). Now (ii) follows from \cite{Sk06, Th. A}. 
The other conclusions hold by Theorem 4.1 since each subalgebra of $H$ is PI.
\endproof

\proclaim
Corollary 4.6.
Suppose that $H$ is a right Noetherian PI Hopf algebra{\rm,} finitely generated 
as an ordinary algebra. Then all conclusions of Theorem\/ {\rm4.5} are true. 
Also{\rm,} each Hopf subalgebra of $H$ is Noetherian.
\endproclaim

\Proof.
It was proved by Anan'in \cite{An92} that any right Noetherian finitely 
generated PI algebra $A$ is residually finite dimensional and, moreover, $A$ 
is \emph{representable} in the sense that $A$ embeds in a finite-dimensional 
algebra over an extension field of the ground field. This property implies the 
ACC on right and left annihilators. Thus $H$ satisfies the assumptions of 
Theorem 4.5.

If $A$ is a Hopf subalgebra of $H$, then, since $H$ is left faithfully $A$-flat, 
the map $I\mapsto IH$ embeds the lattice of right ideals of $A$ into that of 
right ideals of $H$. Hence $A$ is right Noetherian. Since the antipode of $A$ 
is bijective, again by \cite{Sk06, Th. A}, the algebra $A$ is also left 
Noetherian.
\endproof

\proclaim
Lemma 4.7.
Let $H$ be as in Theorem\/ {\rm4.1}. Then the lattice of PI Hopf subalgebras 
of $H$ embeds into the lattice of left ideals of $H$.
\endproclaim

\Proof.
For a PI Hopf subalgebra $A$ of $H$ put $A^+=\{a\in A\mid\ep(a)=0\}$. Then 
$HA^+$ is a left ideal and a coideal of $H$. Since $H$ is faithfully flat over 
$A$ by Theorem 4.1, we can apply Takeuchi's result \cite{Tak79, Th. 1} which 
shows that $A$ can be reconstructed from the quotient coalgebra $H/HA^+$ as
$$
A=\{h\in H\mid(\pi\ot\id)\De(h)=\pi(1)\ot h\ \hbox{ in $(H/HA^+)\ot H$}\}
$$
where $\pi:H\to H/HA^+$ is the canonical map. In other words,
$$
A=\{h\in H\mid\De(h)-1\ot h\in HA^+\ot H\}.
$$
It follows that the assignment $A\mapsto HA^+$ defines an injective 
inclusion-preserving map from the set of PI Hopf subalgebras to the set 
of left ideals of $H$.
\endproof

As an application of the previous results we consider the relationship between 
Noetherianness and finite generation. It was observed by Molnar \cite{Mol75} 
that a commutative Hopf algebra is Noetherian if and only if it is finitely 
generated. Motivated by this fact, Wu and Zhang asked whether every Noetherian 
Hopf algebra is affine, i.e., finitely generated as an ordinary algebra 
\cite{Wu-Zh03, Question 5.1}. Faithful flatness over Hopf subalgebras is an 
essential argument for deriving such a conclusion. We have established this 
property for a certain class of Hopf algebras, but this property alone is not 
sufficient. We are still able to state a somewhat weaker result:

\proclaim
Proposition 4.8.
Let $H$ be a residually finite-dimensional Noetherian PI Hopf algebra.
Then $H$ is finitely generated as a Hopf algebra.
\endproclaim

\Proof.
By \cite{Sk99, Th. 4.5} $H$ has an Artinian classical quotient ring, i.e., $H$ 
satisfies the assumptions of Theorem 4.1. This enables us to use Lemma 4.7. 
Note that all subalgebras of $H$ are PI. Since the ACC is satisfied in the 
lattice of left ideals of $H$, it is also satisfied in the lattice of Hopf 
subalgebras of $H$. It follows that there exists a largest finitely generated 
Hopf subalgebra $A$. Clearly $A$ must contain all elements of $H$, i.e., $A=H$.
\endproof

If in Proposition 4.8 each element of $H$ is contained in a finite-dimensional
subspace stable under the action of the antipode $S$, then $H$ will be finitely 
generated as an ordinary algebra. One particular case where this happens is 
presented below:

\proclaim
Corollary 4.9.
If $H$ is a coquasitriangular residually finite-dimensional Noetherian PI Hopf 
algebra{\rm,} then $H$ is finitely generated as an ordinary algebra.
\endproclaim

\Proof.
It is well-known that $S^2(C)=C$ for any subcoalgebra $C\sbs H$. This is a 
consequence of the fact that the category of right comodules for any 
coquasitriangular Hopf algebra is braided, which implies that $V^{**}\cong V$ 
for each finite-dimensional right $H$-comodule $V$. Moreover, $S^2$ is a 
coinner operator according to Doi \cite{Doi93, Th. 1.3} and Schauenburg 
\cite{Scha92, Lemma 3.3.2}. This means that there exists a convolution 
invertible linear map $\la:H\to k$ such that $S^2(h)=\la\rhu h\lhu\la^{-1}$ 
for all $h\in H$.

In view of Proposition 4.8 we can find a finite-dimensional subcoalgebra $C$ 
which generates $H$ as a Hopf algebra. The subcoalgebra $C'=C+S(C)$ is finite 
dimensional and $S$-invariant. Since the subalgebra of $H$ generated by $C'$ 
is a Hopf subalgebra, it coincides with $H$.
\endproof

\references
\nextref An92
\auth{A.Z.,Anan'in}
\paper{Representability of Noetherian finitely generated algebras}
\journal{Arch. Math.}
\Vol{59}
\Year{1992}
\Pages{1-5}

\nextref Br84
\auth{A.,Braun}
\paper{The nilpotency of the radical in a finitely generated P.I. ring}
\journal{J.~Algebra}
\Vol{89}
\Year{1984}
\Pages{375-396}

\nextref Br07
\auth{K.A.,Brown}
\paper{Noetherian Hopf algebras}
\journal{Turkish J. Math. suppl.}
\Vol{31}
\Year{2007}
\Pages{7-23}

\nextref Chir
\auth{A.,Chirvasitu}
\paper{Cosemisimple Hopf algebras are faithfully flat over Hopf subalgebras}
\journal{Algebra Number Theory}
\Vol{8}
\Year{2014}
\Pages{1179-1199}

\nextref Coh85
\auth{M.,Cohen}
\paper{Hopf algebras acting on semiprime algebras}
\InBook{Group actions on rings}
\BkSer{Contemp. Math.}
\BkVol{43}
\publisher{Amer. Math. Soc.}
\Year{1985}
\Pages{49-61}

\nextref Dem-G
\auth{M.,Demazure;P.,Gabriel}
\book{Groupes Alg\'ebriques I}
\publisher{Masson}
\Year{1970}

\nextref Doi93
\auth{Y.,Doi}
\paper{Braided bialgebras and quadratic bialgebras}
\journal{Comm. Algebra}
\Vol{21}
\Year{1993}
\Pages{1731-1749}

\nextref Goo
\auth{K.R.,Goodearl}
\book{Nonsingular Rings and Modules}
\publisher{Marcel Dekker}
\Year{1976}

\nextref Goo-W
\auth{K.R.,Goodearl;R.B.,Warfield Jr.}
\book{An Introduction to Noncommutative Noetherian Rings}
Second edition,
\publisher{Cambridge Univ. Press}
\Year{2004}

\nextref Gro
\auth{F.D.,Grosshans}
\book{Algebraic Homogeneous Spaces and Invariant Theory}
\BkSer{Lecture Notes Math.}
\BkVol{1673}
\publisher{Springer}
\Year{1997}

\nextref Her-S64
\auth{I.N.,Herstein;L.,Small}
\paper{Nil rings satisfying certain chain conditions}
\journal{Canad. J. Math.}
\Vol{16}
\Year{1964}
\Pages{771-776}

\nextref Kr-L
\auth{G.R.,Krause;T.H.,Lenagan}
\book{Growth of Algebras and Gelfand-Kirillov Dimension}
Revised edition,
\publisher{Amer. Math. Soc.}
\Year{2000}

\nextref Mal43
\auth{A.I.,Malcev}
\paper{On representations of infinite algebras\inRus}
\journal{Mat. Sbornik}
\Vol{13}
\Year{1943}
\Pages{263-286}

\nextref Ma91
\auth{A.,Masuoka}
\paper{On Hopf algebras with cocommutative coradicals}
\journal{J.~Algebra}
\Vol{144}
\Year{1991}
\Pages{451-466}

\nextref Ma-W94
\auth{A.,Masuoka;D.,Wigner}
\paper{Faithful flatness of Hopf algebras}
\journal{J.~Algebra}
\Vol{170}
\Year{1994}
\Pages{156-164}

\nextref Mc-R
\auth{J.C.,McConnell;J.C.,Robson}
\book{Noncommutative Noetherian Rings}
Revised edition,
\publisher{Amer. Math. Soc.}
\Year{2001}

\nextref Mol75
\auth{R.K.,Molnar}
\paper{A commutative Noetherian Hopf algebra over a field is finitely generated}
\journal{Proc. Amer. Math. Soc.}
\Vol{51}
\Year{1975}
\Pages{501-502}

\nextref Mo
\auth{S.,Montgomery}
\book{Hopf Algebras and their Actions on Rings}
\publisher{Amer. Math. Soc.}
\Year{1993}

\nextref Mo-Sch95
\auth{S.,Montgomery;H.-J.,Schneider}
\paper{Hopf crossed products, rings of quotients, and prime ideals}
\journal{Adv. Math.}
\Vol{112}
\Year{1995}
\Pages{1-55}

\nextref Mu-Sch99
\auth{E.F.,M\"uller;H.-J.,Schneider}
\paper{Quantum homogeneous spaces with faithfully flat module structures}
\journal{Israel J. Math.}
\Vol{111}
\Year{1999}
\Pages{157-190}

\nextref Nich-Z89
\auth{W.D.,Nichols;M.B.,Zoeller}
\paper{A Hopf algebra freeness theorem}
\journal{Amer. J. Math.}
\Vol{111}
\Year{1989}
\Pages{381-385}

\nextref Pro
\auth{C.,Procesi}
\book{Rings with Polynomial Identities}
\publisher{Marcel Dekker}
\Year{1973}

\nextref Rad77
\auth{D.E.,Radford}
\paper{Pointed Hopf algebras are free over Hopf subalgebras}
\journal{J.~Algebra}
\Vol{45}
\Year{1977}
\Pages{266-273}

\nextref Row
\auth{L.H.,Rowen}
\book{Ring Theory}
\publisher{Academic Press}
\Year{1988}

\nextref Scha92
\auth{P.,Schauenburg}
\book{On Coquasitriangular Hopf algebras and the Quantum Yang-Baxter Equation}
\BkSer{Algebra Berichte 67}
\publisher{Reinhard Fisher}
\Year{1992}

\nextref Scha00
\auth{P.,Schauenburg}
\paper{Faithful flatness over Hopf subalgebras: counterexamples}
\InBook{Interactions between ring theory and representations of algebras}
\publisher{Marcel Dekker}
\Year{2000}
\Pages{331-344}

\nextref Schn92
\auth{H.-J.,Schneider}
\paper{Normal basis and transitivity of crossed products for Hopf algebras}
\journal{J.~Algebra}
\Vol{152}
\Year{1992}
\Pages{289-312}

\nextref Schn93
\auth{H.-J.,Schneider}
\paper{Some remarks on exact sequences of quantum groups}
\journal{Comm. Algebra}
\Vol{21}
\Year{1993}
\Pages{3337-3357}

\nextref Sk06
\auth{S.,Skryabin}
\paper{New results on the bijectivity of antipode of a Hopf algebra}
\journal{J.~Algebra}
\Vol{306}
\Year{2006}
\Pages{622-633}

\nextref Sk07
\auth{S.,Skryabin}
\paper{Projectivity and freeness over comodule algebras}
\journal{Trans. Amer. Math. Soc.}
\Vol{359}
\Year{2007}
\Pages{2597-2623}

\nextref Sk08
\auth{S.,Skryabin}
\paper{Projectivity of Hopf algebras over subalgebras with semilocal central localizations}
\journal{J.~$K$-Theory}
\Vol{2}
\Year{2008}
\Pages{1-40}

\nextref Sk10
\auth{S.,Skryabin}
\paper{Models of quasiprojective homogeneous spaces for Hopf algebras}
\journal{J.~Reine Angew. Math.}
\Vol{643}
\Year{2010}
\Pages{201-236}

\nextref Sk99
\auth{S.,Skryabin}
\paper{Flatness of Noetherian Hopf algebras over coideal subalgebras}
arXiv:2001. 02848,
to appear in \journal{Algebras Represent. Theory}.

\nextref Sk-Oy06
\auth{S.,Skryabin;F.,Van Oystaeyen}
\paper{The Goldie theorem for $H$-semiprime algebras}
\journal{J.~Algebra}
\Vol{305}
\Year{2006}
\Pages{292-320}

\nextref St
\auth{B.,Stenstr\"om}
\book{Rings of Quotients}
\publisher{Springer}
\Year{1975}

\nextref Tak79
\auth{M.,Takeuchi}
\paper{Relative Hopf modules---equivalences and freeness criteria}
\journal{J.~Algebra}
\Vol{60}
\Year{1979}
\Pages{452-471}

\nextref Wu-Zh03
\auth{Q.-S.,Wu;J.J.,Zhang}
\paper{Noetherian PI Hopf algebras are Gorenstein}
\journal{Trans. Amer. Math. Soc.}
\Vol{355}
\Year{2003}
\Pages{1043-1066}

\endreferences
\bye